\documentclass{birkjour_t2}

\usepackage{amsthm}
\newtheorem{theorem}{Theorem}[section]
\newtheorem{lemma}[theorem]{Lemma}
\theoremstyle{definition}
\newtheorem{definition}[theorem]{Definition}
\theoremstyle{remark}
\newtheorem{remark}[theorem]{Remark}
\renewcommand{\div}{ {\rm div }  }

\def\u{\mathbf{u}}
\def\f{\frac}
\def\x{\mathbf{x}}
\def\d{\mathrm{d}}
\def\pf{\it{Proof.}\rm\quad}

\begin{document}

\title[Free boundary value problem to
Navier-Stokes-Poisson equations ]
 {Free boundary value problem to 3D spherically symmetric\\ compressible
Navier-Stokes-Poisson equations}

\author[H.-H. Kong]{Huihui Kong$^{1,2}$}

\email{konghuihuiking@126.com }
\author[H.-L. Li]{Hai-Liang Li$^{1}$}
\email{hailiang.li.math@gmail.com }

 \address{%
 \small $^1$School of Mathematical Sciences,\\
 Capital Normal University,
 Beijing 100048, P.R.China\\
\small $^2$ Beijing Computational Science Research Center, \\
Beijing 100193, P. R. China}




\begin{abstract}
In the paper, we consider the free boundary value problem to 3D spherically symmetric compressible
isentropic Navier-Stokes-Poisson equations for self-gravitating gaseous stars with $\gamma$-law pressure
density function for $\frac65 <\gamma \leq \frac43$. For stress free boundary condition and zero flow density
 continuously  across the free boundary, the global existence of spherically symmetric weak solutions is shown, and the regularity
and long time behavior of global solution are investigated for spherically symmetric initial data with the total mass smaller than a critical mass.
\end{abstract}

\maketitle
\section{Introduction and Main Results}
The motion of self-gravitating viscous gaseous stars can be described by the
compressible Navier-Stokes-Poisson (NSP) system  in $\mathbb{R}^3$ :
\begin{equation}\label{1.1}
\left\{%
\begin{array}{l}
  \rho_t+\div (\rho\mathbf{u})=0,
  \\
  (\rho\mathbf{u})_t+\div (\rho\mathbf{u}\otimes\mathbf{u})
  +\nabla P(\rho)=\mu\triangle\mathbf{u}+(\lambda+\mu)\nabla\div\u + \rho \nabla \Phi,\\
-\triangle\Phi=4\pi\rho
\end{array}%
\right.
\end{equation}
where $(\x,t)\in\mathbb{R}^3\times \mathbb{R}_{+}, $ $\rho(\x,t),\u(\x,t)=(u_1(\x,t),u_2(\x,t),u_3(\x,t)),~P(\rho)$ and $\Phi$ denote the density, velocity, pressure and the potential function of the self-gravitational force respectively. The constants $\mu$ and $\lambda$
are Lam\'{e} viscosity coefficients satisfying
\begin{equation}
   \mu>0,\ 2\mu+3\lambda\geq 0.
\end{equation}
We assume the polytropic gas pressure
\begin{equation}
    P(\rho)=\kappa \rho^\gamma,\ \gamma>1,
\end{equation}
where $\gamma>1$ is the adiabatic exponent and $\kappa>0$ is an entropy constant which is to be one for simplicity. In the motion of gaseous stars,
the different value of $\gamma$ may affect the existence, uniqueness and stability of stationary solutions.
 For the spherically symmetric motion, stationary solution $(\bar{\rho}(r),0)$ with non-moving gaseous sphere, satisfies the following:
\begin{equation}\label{s1}
    \partial_r\bar{\rho}^\gamma=-\f{4\pi\bar{\rho}}{r^2}\int_{0}^{r}\bar{\rho}s^2\d s.
\end{equation}
The ordinary differential equation \eqref{s1} can be transformed into the famous Lane-Emden equation \cite{C1938}.
 For given finite total mass, there exists at least one compactly supported stationary solution to \eqref{s1} for $\f{6}{5}<\gamma<2$,
 refer to \cite{L1997}. For $\f{4}{3}<\gamma<2$, every stationary solution is compactly supported  and unique. For $\gamma=\f{6}{5}$, there is a unique solution with infinite support, and it can be written explicitly in terms of the Lane-Emden function. On the other hand, for $1<\gamma<\f{6}{5}$, there are no stationary solutions with finite total mass.
Recently, many important study on the asymptotic stability/instability of stationary solutions has been made,
for instance, in \cite{DLYY,L1997,JT2013,JJ2008,LXZ,LXZ1,ZF2009}. Due to the balance between
flow pressure and gravity force, it was shown that there is a critical value $\gamma_c=\f{4}{3}$
in spatial three-dimension, the stationary Lane-Emden solution is expected to be stable for  $\gamma>\f{4}{3}$ and
 instable for  $\gamma<\f{4}{3}$. Indeed, Jang and Tice \cite{JT2013} prove the instability theory of the NSP equations
for $\f{6}{5}<\gamma<\f{4}{3}$, and the nonlinear asymptotic stability of the Lane-Emden solutions for
the viscous gaseous star is established by Luo-Xin-Zeng \cite{LXZ,LXZ1} for $\f{4}{3}<\gamma<2$.
For $\gamma=\f{6}{5}$, the nonlinear instability is proved for gravitational Euler-Poisson system in~\cite{JJ2008}.
\par

The free boundary value problem (FBVP) for the compressible NSP equations  which involves the influence of the vacuum state on the existence and dynamics of solutions has attracted lots of research interests and been studied extensively, refer to \cite{SG,JJ2010,JT2013} and the references therein. For  instance, Jang\cite{JJ2010} establishes the local-in-time well-posedness of spherically symmetric strong solutions to the free boundary value problem for $\f{6}{5}<\gamma<2$.
The global spherically symmetric weak solutions are constructed in \cite{SG} for $\gamma>\f{4}{3}$
where the choice of $\gamma>\f{4}{3}$ plays a critical role on the energy estimates to ensure that
the negative gravitational energy can be dominated by the positive kinetic-internal/dissipation energy.
The global existence of a spherically symmetric entropy weak solution for the compressible NSP system
with density-dependent viscosity coefficients to the FBVP is shown for general initial data with finite entropy
 when the density changes discontinuously across the interfaces separating the fluid
 and vacuum in \cite{DL2015} for $\f{6}{5}<\gamma\leq\f{4}{3}$. However, it is not known
 yet whether there may exist or not any solution globally in time to NSP~\eqref{1.1} for $\f{6}{5}<\gamma\leq\f{4}{3}$ with
 the stress-free boundary condition and the zero flow density continuously across the free boundary. And it is also interesting to investigate the long time behaviors of global solution to the FBVP.\par
In this paper, we investigate the FBVP for the compressible NSP system (\ref{1.1}) for viscous gaseous stars with the stress-free boundary condition and the zero flow density across the free boundary. For spherically symmetric initial data with finite mass (smaller than a critical mass) and energy, we prove the global existence of spherically symmetric weak solutions to the FBVP problem for (\ref{1.1}), establish the regularities of solution and the positivity of flow density, and obtain the expanding rate of the domain occupied by the fluid for $\f{6}{5}<\gamma\leq \f{4}{3}$.

To state the main results, let us consider the spherically symmetric solution $(\rho,\mathbf{u})$ to (\ref{1.1}) in $\mathbb{R}^3$ so that $$\rho(\mathbf{x},t)=\rho(r,t),\ \mathbf{u}=u(r,t)\frac{\mathbf{x}}{r},\ r=|\mathbf{x}|,\ $$
 and (\ref{1.1}) are changed to
 \begin{eqnarray}\label{1.2}
 \left\{%
 \begin{array}{l}
   \rho_t+(\rho u)_r+\frac{2\rho u}{r}=0,\\
   (\rho u)_t+(\rho u^2+\rho^\gamma)_r+\frac{2\rho u^2}{r}=(\lambda+2\mu)(u_r+\frac{2u}{r})_r-\f{4\pi\rho}{r^2}\int_{0}^{r}\rho s^2\d s,
 \end{array}%
 \right.
 \end{eqnarray}
for $(r,t)\in\Omega_T$ with
 $$\Omega_T=\{(r,t)|0\leq r\leq a(t),\ 0\leq t\leq T\}.$$
 The initial data is taken as
\begin{equation}\label{1.3}
    (\rho,\rho u)(r,0)=(\rho_0,m_0)(r):=(\rho_0,\rho_0 u_0)(r),\ r\in(0, a_0).
\end{equation}
At the center of symmetry we impose the Dirichlet boundary condition
\begin{equation}\label{1.4}
     u(0,t)=0,
\end{equation}
and across the free surface $\partial\Omega_t$ which moves in the radial direction along the particle path $r=a(t)$, the vacuum state appears and the stress-free boundary condition holds
\begin{equation}\label{1.5}
 F(a(t),t)=0,\quad \rho(a(t),t)=0,\quad  t\geq0,
\end{equation}
where $a^\prime(t)=u(a(t),t),\ t>0$\ , $a(0)=a_0>0$ and the stress (effective viscous flux) $F$ is defined by
\begin{equation}\label{Flux}
F=:\rho^\gamma-(\lambda+2\mu) \div\u=\rho^\gamma-(\lambda+2\mu) u_r -(\lambda+2\mu)\frac{2u}{r}.
\end{equation}\par
\begin{definition}\label{def}
  $(\rho,\mathbf{u},a)$ with $\rho\geq0$ a.e. is said to be a weak solution to the free boundary value problem \eqref{1.1}, \eqref{1.2}-\eqref{1.5} on $\Omega_t\times[0,T]$, provided that it holds that
\begin{equation*}
    \rho \in L^\infty(0,T;L^1(\Omega_t)\cap L^\gamma(\Omega_t)),\ \sqrt{\rho}\mathbf{u}\in L^\infty(0,T;L^2(\Omega_t)),
\end{equation*}
\begin{equation}
  \nabla\mathbf{u}\in L^2(0,T;L^2(\Omega_t)),\
a(t)\in H^1([0,T]), 
\end{equation}
and the equations are satisfied in the sense of distribution. Namely, it holds for any $t_2>t_1\geq 0$ and any $\phi\in C^1([0,T]\times\bar{\Omega}_t)$ that
\begin{equation}
\int_{\Omega_t}\rho\phi\mathrm{d}\mathbf{x}|_{t_1}^{t_2}
=\int_{t_1}^{t_2}\int_{\Omega_t}(\rho\phi_t+\rho \mathbf{u}\cdot\nabla\phi)\mathrm{d}\mathbf{x}\mathrm{d}t,
\end{equation}
and for $\psi=(\psi_1,\psi_2,\psi_3)\in C^1([0,T]\times\bar{\Omega}_t)$ satisfying $\psi(\mathbf{x},T)=0$ and $\psi(\mathbf{x},t)=0$ on $\partial\Omega_t$ that
\begin{equation*}
 \int_{\Omega_t}\mathbf{m_0}\cdot\psi(\mathbf{x},0)\d\mathbf{x}
+\int_{0}^{T}\int_{\Omega_t}[\rho\mathbf{u}\cdot\partial_t\psi
+\rho\mathbf{u}\otimes\mathbf{u}:\nabla\psi]\mathrm{d}\mathbf{x}\mathrm{d}t
\end{equation*}
\begin{equation}
 +\int_{0}^{T}\int_{\Omega_t}\rho^\gamma \div \psi\mathrm{d}\mathbf{x}\mathrm{d}t +\int_{0}^{T}\int_{\Omega_t}(\mu \nabla\u:\nabla\psi+(\lambda +\mu) \div\u\div \psi)\mathrm{d}\mathbf{x}\mathrm{d}t=\int_{0}^{T}\int_{\Omega_t}\rho \nabla \Phi ~\psi \mathrm{d}\mathbf{x}\mathrm{d}t,
\end{equation}
where $\Omega_t=\{r\in \mathbb{R_+}|0\leq r\leq a(t)\},$
and
\begin{equation*}
-\triangle\Phi=4\pi\rho~~a.e.
\end{equation*}
where $\Phi$ is defined in $\mathbb{R}^3$ through the Poisson equation with $\Phi\rightarrow 0$ as $|\x|\rightarrow\infty$ and $\rho=0$ in $\mathbb{R}^3 \backslash\Omega_t $ at time $t$.
The free boundary condition \eqref{1.5} is satisfied in the sense of trace.
\end{definition}

Then, we have the main results on global existence and long time behavior of solution to free boundary value problem of \eqref{1.1} for $\f{6}{5}<\gamma\leq\f{4}{3}$ below.
\begin{theorem}\label{thm1}
$\mathrm{(Global\  existence})$  Let $T>0$ and $\f{6}{5}<\gamma\leq\f{4}{3}$.
  Assume that  the spherically symmetric initial data \eqref{1.3} satisfies
the regularity and compatibility conditions
 \begin{equation}\label{t2.9}
    0\leq\rho_0\in L^{1}(\Omega_0)\cap L^{\infty}(\Omega_0),\ (\rho_0)^{k}\in H^{1}(\Omega_0),\ \u_0\in H^1(\Omega_0),
 \end{equation}
 \begin{equation}\label{t2.10}
   \ \rho_0(r)>0,\  r\in(0,a_0),\ \rho_0(a_0)=0,\  u_{0 r }(a_0)+\frac{2u_0(a_0)}{a_0}=0,
 \end{equation}
where $k$ is the constant satisfying $0<k\leq\gamma-\f{1}{2}$. If the mass $M=4\pi\int_{0}^{a_0}\rho_0(r) r^2\mathrm{d}r$ satisfies $M<M_c$
with
 $$ M_c=
\left\{ \aligned
  &(\f{3}{B})^{\f{3}{2}},&\gamma=\f{4}{3} ;\\
  &[\f{4-3\gamma}{\gamma-1}(\f{B}{3})^{-\f{3(\gamma-1)}{4-3\gamma}}]^{\f{4-3\gamma}{5\gamma-6}}(E_0)^{-\f{4-3\gamma}{5\gamma-6}} ,&\gamma\in(\f{6}{5},\f{4}{3}),
 \endaligned\right. $$
and B is a positive constant related to $\gamma$ defined in Lemma \ref{lem2.1},
 then there exists a global spherically symmetric weak solution $$(\rho,\u,a)(\x,t)=(\rho(r,t),u(r,t)\f{|\x|}{r},a(t)),\ r=|\x|,$$
to $(\ref{1.2})$  for $t\in[0,T]$, which satisfies $\rho(r,t)\ge 0$ a.e.  and
\begin{align}\label{enegy1}
   & c_0\leq a(t)\leq C_T,~~ a(t)\in H^1(0,T),\\
   & \int_{0}^{a(t)}(\frac{1}{2}\rho u^2+C_\gamma\rho^\gamma) r^2\mathrm{d}r+ (\lambda+2\mu)\int_{0}^{t}\int_{0}^{a(\tau)}(u_r+\frac{2u}{r})^2 r^2\mathrm{d}r\d\tau\leq E_0,
\end{align}
where $c_0>0$, $C_T>0$ are two constants, $E_0:=\int_{0}^{a_0}(\frac{1}{2}\rho_0 u_0^2+\frac{1}{\gamma-1}\rho_0^\gamma) r^2\mathrm{d}r$ and $C_\gamma=\f{4-3\gamma}{\gamma-1},\gamma\in(\f{6}{5},\f{4}{3}); ~ C_\gamma =3-BM^{\f{2}{3}}>0$  for $\gamma=\f{4}{3}$.\par
Furthermore, the solution $(\rho, u ,a)$ satisfies the following properties:

$\mathit{(i)}\ (\mbox{Transport Property})$ For any $r_i\in (0,a_0]$, there exist positive constants $C_{x_i,T},~c_{x_i,T}$ and $~c$ such that
 \begin{gather}\label{2.13a}
\rho_0(r_i)e^{-c_{x_i,T}/(\lambda+2\mu)}\leq \rho (r_{x_i}(t),t)\leq \rho_0(r_i)e^{C_{x_i,T}/(\lambda+2\mu)},\  \ t\in[0,T],
\\
 c x_0^{\frac{\gamma}{3(\gamma-1)}}\leq r_{x_0}(t)\leq a(t),\ t\in[0,T],
\\
 c(x_2-x_1)^{\frac{\gamma}{\gamma-1}}\leq r_{x_2}^3(t)- r_{x_1}^3(t),\ t\in[0,T],
\end{gather}
where $r_{x_i}(t), \ i=0,1,2,$ is the particle path defined by
$ \f{\d r_{x_i}(t)}{\d t}=u(r_{x_i}(t),t)$ with $r_{x_i}(0)=r_i\in(0,a_0]$ and $x_i=\f{M}{4\pi}-\int_{r_i}^{a_0}\rho_0 r^2\mathrm{d}r,$   $C_{x_i,T},~c_{x_i,T}\rightarrow+\infty$ as $x_i\rightarrow0$.

$\mathrm{(ii)\ ( Interior\  regularity)}$ If the initial velocity also satisfies $u_0\in H^2([r_0^-,r_b^+])$ for any $0<r_0^-<r_0<r_b<r_b^+\leq a_0$.
Then, the following interior regularities hold
 \begin{equation}\label{t2.14}
\left\{
  \begin{aligned}
    &(\rho, u)\in C([r_{x_0}(t),r_{x_b}(t)]\times[0,T]), &\hbox{}
    \\
    & \rho\in L^\infty(0,T; H^1([r_{x_0}(t),r_{x_b}(t)])), u \in L^\infty(0,T; H^2([r_{x_0}(t),r_{x_b}(t)])),&\hbox{}
    \\
    & \rho_t\in L^\infty(0,T; L^2([r_{x_0}(t),r_{x_b}(t)]))\cap L^2(0,T; H^1([r_{x_0}(t),r_{x_b}(t)])),&\hbox{}
    \\
    & u_t\in L^\infty(0,T; L^2([r_{x_0}(t),r_{x_b}(t)]))\cap L^2(0,T; H^1([r_{x_0}(t),r_{x_b}(t)])) & \hbox{}
  \end{aligned}
\right.
\end{equation}
where $r_{x_0}(t)$ is the particle path defined as above and $r_{x_b}(t)$ is the particle path with $r_{x_b}(0)=r_b$ and $x_b=\f{M}{4\pi}-\int_{r_b}^{a_0}\rho_0 r^2\mathrm{d}r.$

$\mathrm{(iii)\ ( Boundary\  regularity)}$
It holds near the free boundary $r=a(t)$ that
\begin{align}\label{2.14a}
 \nonumber  & \|(\rho^k,u)(t)\|_{H^1(\Omega_\eta)}+\|F(t)\|_{L^2(\Omega_\eta)}+\|\sqrt{\rho}\dot{u}\|_{L^2(0,T; L^2(\Omega_\eta))}\\
&+\|F\|_{L^2(0,T; H^{1}(\Omega_\eta))}+\|u\|_{L^2(0,T;H^{2}(\Omega_\eta))}+\|a\|_{H^1([0,T])}\leq C_{T}\delta_0,
\end{align}
with $\delta_0=:\|\rho_0\|_{L^{\infty}([0,a_0])}+\|u_0\|_{H^1([0,a_0])}+\|\rho^k_0\|_{H^{1}([0,a_0])}$ and $\Omega_\eta=(a(t)-\eta,a(t))$ for some small constant $\eta>0$.
In addition, if the initial data $(\rho_0, u_0)$ satisfy $u_0\in H^2([a_0-\eta,a_0]),\ \rho_0^{-\f{1}{2}}\partial_r^2u_0\in L^2([a_0-\eta,a_0])$ and compatibility condition, then
\begin{align}\label{2.15a}
 &  \|\sqrt{\rho}\dot{u}(t)\|_{L^{2}(\Omega_\eta)} +\|u(t)\|_{H^{2}(\Omega_\eta)}+\|\rho^{-\f{1}{2}}\partial_r^2u\|_{L^2(\Omega_\eta)}
+\|F(t)\|_{H^{1}(\Omega_\eta)} + \|a\|_{H^2([0,T])}\leq C_{T}\delta_1,
\end{align}
with $\dot{u}=u_t+uu_r$ and $\delta_1=\|\rho_0\|_{L^{\infty}([0,a_0])}+\|u_0\|_{H^1([0,a_0])}+\|\rho_0^{-\f{1}{2}}\partial_r^2u_0\|_{L^2([a_0-\eta,a_0])}+\|\rho^k_0\|_{H^{1}([0,a_0])}.$

\end{theorem}
\begin{theorem}\label{thm2}
Let $T>0,\ \f{6}{5}<\gamma\leq\f{4}{3}$ and $(\rho,u,a)$ be any global  (strong or weak) solution to the FBVP (\ref{1.2}) for $t\in[0,T]$ with $F=\rho^\gamma-(\lambda+2\mu) u_r-(\lambda+2\mu)\frac{2 u}{r}\in$
 $ L^2(0,T;H^1(\Omega_\eta))$ and $\Omega_\eta=(a(t)-\eta,a(t))$ for some small constant $\eta>0.$ If $M<\overline{M}<M_c$ with
 $$ \overline{M}=
\left\{ \aligned
  &(\f{3}{2B})^{\f{3}{2}},&\gamma=\f{4}{3} ;\\
  &[\f{4-3\gamma}{\gamma-1}(\f{B}{3})^{-\f{3(\gamma-1)}{4-3\gamma}}]^{\f{4-3\gamma}{5\gamma-6}}(lE_0)^{-\f{4-3\gamma}{5\gamma-6}} ,&\gamma\in(\f{6}{5},\f{4}{3}),
 \endaligned\right. $$
 for some $l>1$, then
 \begin{equation}\label{2.140}
\int_{0}^{a(t)}(\frac{1}{2}\rho u^2+\f{1}{2(\gamma-1)}\rho^\gamma) r^2\mathrm{d}r+ (\lambda+2\mu)\int_{0}^{t}\int_{0}^{a(\tau)}(u_r+\frac{2u}{r})^2 r^2\mathrm{d}r\d\tau\leq E_0,
\end{equation} and then for any $t>0$,
\begin{equation}\label{2.141}
 \frac{1}{a^3_1(t)}\int_{0}^{a(t)}\rho^{\gamma} r^{2}\mathrm{d}r\leq C(1+t)^{7-6\gamma},
\end{equation}
therefore,
\begin{equation}\label{2.14}
a_1(t)=\max_{s\in[0,t]}a(s)\geq C(1+t)^{\frac{6\gamma-7}{3\gamma}},~~\f{6}{5}<\gamma\leq\f{4}{3},
\end{equation}
where C is a constant independent of time.\par
In particular, it holds for $\gamma=\f{4}{3}$ that
\begin{equation}\label{2.151}
  \f{1}{a^3(t)}\int_{0}^{a(t)}\rho^{\gamma}r^2 \d r\leq C (1+t)^{-1},
\end{equation}
therefore,
\begin{equation}\label{2.15}
a(t)\geq C(1+t)^{\f{1}{4}},
\end{equation}
and for $\gamma\in (\f{6}{5},\f{4}{3})$ and any $\beta\in (\f{2(4-3\gamma)}{3},\f{1}{3\gamma})$, there exists time sequence $\{t_n\}$ such that
\begin{equation}\label{ja14}
   \f{a(t_n)}{(1+t_n)^\beta}\rightarrow +\infty,~~t_n\rightarrow +\infty.
\end{equation}
\end{theorem}
{\remark As shown in \cite{L1997}, for given finite total mass $M$, there exists at least one compactly supported stationary solution $(\bar{\rho}(r),0)$ satisfying \eqref{s1}
in a finite domain $[0,\bar{a}]$ for $\f{6}{5}<\gamma\leq\f{4}{3}$. And it is not difficult to analyze that the steady state $(\bar{\rho}(r),0)$ does
 not fit the phenomenon shown in Theorem \ref{thm2}.
Indeed, it is easy to obtain that the energy of $(\bar{\rho}(r),0)$ satisfies $\frac{1}{\gamma-1}\int_{0}^{\bar{a}}\bar{\rho}^\gamma r^2\mathrm{d}r
-4\pi\int_{0}^{\bar{a}}\bar{\rho} r\int_{0}^{r}\bar{\rho}s^2\d s\mathrm{d}r=\frac{4-3\gamma}{\gamma-1}\int_{0}^{\bar{a}}\bar{\rho}^\gamma r^2\mathrm{d}r
< \frac{1}{2(\gamma-1)}\int_{0}^{\bar{a}}\bar{\rho}^\gamma r^2\mathrm{d}r$
for $\f{6}{5}<\gamma\leq\f{4}{3}$. But for any solution which satisfies the assumption of Theorem \ref{thm2}, its energy satisfies
$\frac{1}{\gamma-1}\int_{0}^{a(t)}\rho^\gamma r^2\mathrm{d}r
-4\pi\int_{0}^{a(t)}\rho r\int_{0}^{r}\rho s^2\d s\mathrm{d}r$ $\geq$ $\f{1}{2(\gamma-1)}\int_{0}^{a(t)}\rho^\gamma r^2\mathrm{d}r$
(which will be proved in Lemma \ref{lem2.1}). This is a contradiction provided that $(\bar{\rho}(r),0)$ satisfies the assumption of Theorem \ref{thm2}.
In addition, one can conclude from $\eqref{2.141}$ that the mean value of total pressure is dispersive for any dynamical solution which
is different from the stationary solution.\\}

The rest part of the paper is arranged as follows. In Sect. 2, the uniform a-priori estimates
 of global approximate solutions are established and the Theorem \ref{thm1} on global existence
 of spherically symmetric solution to original problem is shown. Particularly, we make basic energy estimate
, the integrability of the pressure and the bounds of density in the Eulerian coordinates. After that,
 by the Lagrangian coordinates transform we translates the moving domain into a fixed domain to estimate the higher
 regularity of the approximate solutions near the free boundary and in the interior domain.  In Sect. 3, the Theorem \ref{thm2} on the long time expanding rate of the domain is established.

\section{A-Priori Estimates}
\label{sect2}
To prove Theorem~\ref{thm1}, let us consider the following modified free boundary value problem:
 \begin{eqnarray}\label{1.2a}
 \left\{%
 \begin{array}{l}
   \rho_t+(\rho u)_r+\frac{2\rho u}{r}=0,\\
   (\rho u)_t+(\rho u^2+\rho^\gamma)_r+\frac{2\rho u^2}{r}=(\lambda+2\mu)(u_r+\frac{2u}{r})_r-\f{4\pi\rho}{r^2}\int_{\varepsilon}^{r}\rho s^2\d s,
 \end{array}%
 \right.
 \end{eqnarray}
for $(r,t)\in\Omega^\varepsilon_T=\{(r,t)|\varepsilon\leq r\leq a(t),\ 0\leq t\leq T\}$ with the following initial data and boundary condition for any fixed small $\varepsilon>0$
\begin{equation}\label{3.1}
   (\rho^\varepsilon,u^\varepsilon)(r,0)=(\rho^\varepsilon_0,u^\varepsilon_0)(r),\quad \varepsilon\leq r\leq a_0,
\end{equation}
\begin{equation}\label{3.2}
    u^\varepsilon(\varepsilon,t)=0,\ ((\rho^\varepsilon)^\gamma-(\lambda+2\mu) u^\varepsilon_r-(\lambda+2\mu) \frac{2u^\varepsilon}{r})(a^\varepsilon(t),t)=0,\quad t>0,
\end{equation}
where $(a^{\varepsilon})^\prime(t)=u^\varepsilon((a^\varepsilon(t),t),\ t>0$ and $a^\varepsilon(0)=a_0$. We should establish the global existence and the uniformly a priori estimates of the approximate solution sequence with respect to $\varepsilon>0$. Without the loss of generality, one can assume that the initial data is smooth enough and consistent with the boundary value (\ref{3.2}) to the higher order. The uniform estimates of the approximate solutions will be made. From simplicity, we omit the subscript $\varepsilon $ below.

To deal with the motion of free boundary, it is  convenient to investigate the approximate FBVP problem (\ref{1.2a})-(\ref{3.2}) in Lagrangian coordinates . Let $(\rho,u,a)$ be any strong solution to the FBVP \eqref{1.2a}-\eqref{3.2}. By the conservation of mass $$M=4\pi\int_{\varepsilon}^{a(t)}\rho r^2\mathrm{d}r=4\pi\int_{\varepsilon}^{a_0}\rho_0r^2\mathrm{d}r,$$ define the Lagrangian coordinates transform
\begin{equation}\label{3.7}
    x(r,t)=\int_{\varepsilon}^{r}\rho y^2\mathrm{d}y,\ \tau=t,
 \end{equation}
for $(r,t)\in\Omega^\varepsilon_T=\{(r,t)|\varepsilon\leq r\leq a(t),\ 0\leq t\leq T\}$ which translates the domain $\Omega^\varepsilon_T$ into $[0,\f{M}{4\pi}]\times[0,T]$ and satisfies
\begin{equation}\label{3.8}
    \frac{\partial x}{\partial r}=\rho r^2,\ \frac{\partial x}{\partial t}=-\rho u r^2,\  \frac{\partial \tau}{\partial r}=0,\  \frac{\partial \tau}{\partial t}=1,
 \end{equation}
 and
 \begin{equation}\label{3.9}
    r^3(x,\tau)=\varepsilon^3+3\int_{0}^{x}\frac{1}{\rho}(y,\tau)\d y=a^3(t)-3\int_{x }^{\f{M}{4\pi}}\frac{1}{\rho}(y,\tau)\d y,\ \frac{\partial r}{\partial \tau}=u.
 \end{equation}
The free boundary problem (\ref{1.2a})-(\ref{3.2}) is changed to
\begin{eqnarray}\label{3.10}
  \left\{%
  \begin{array}{l}
  \rho_{\tau}+\rho^2(u r^2)_x=0,
  \\
  u_{\tau}+r^2(\rho^\gamma-(\lambda+2\mu)\rho(u r^2)_x)_x+\f{4\pi x}{r^2}=0,
\end{array}%
\right.
\end{eqnarray}
The initial data and boundary conditions are given by
\begin{eqnarray}\label{3.11}
  \nonumber &&(\rho,u)(x,0)= (\rho_0,u_0)(x) \\
 &&u(0,\tau)= 0,\ (\rho^\gamma-(\lambda+2\mu)\rho r^2u_x-(\lambda+2\mu)\frac{2u}{r})(\f{M}{4\pi},\tau)=0,
\end{eqnarray}
where $r=r(x,\tau)$ is defined by
\begin{equation}\label{3.12}
   \frac{\partial}{\partial\tau}r(x,\tau)=u(x,\tau),\ x\in[0,\f{M}{4\pi}],\ \tau\in[0,T],
\end{equation}
 and the boundary $x=\f{M}{4\pi}$ corresponds to the free boundary $a(\tau)=r(\f{M}{4\pi},\tau)$ determined by
\begin{equation}\label{3.13}
    \frac{\mathrm{d}}{\mathrm{d}\tau}a(\tau)=u(\f{M}{4\pi},\tau),\ \tau\in[0,T] ;\ a(0)=a_0.
\end{equation}
\par
In the rest part of this section, we make basic energy estimate, the integrability of the pressure and the bounds of density in the Eulerian coordinates and make boundary and interior estimates in Lagrangian coordinates.
\subsection{Estimates in Eulerian coordinates}
\begin{lemma}\label{lem2.1}(Basic energy estimate)
Let $T>0$, $\f{6}{5}<\gamma\leq\f{4}{3}$ and $(\rho,u,a)$  be any strong solution to the FBVP \eqref{1.2a}-\eqref{3.2}. If the mass $M<M_c$, then it holds that
\begin{equation}\label{3.19}
\int_{\varepsilon}^{a(t)}(\frac{1}{2}\rho u^2+C_\gamma\rho^\gamma) r^2\mathrm{d}r+ (\lambda+2\mu)\int_{0}^{t}\int_{\varepsilon}^{a(\tau)}(u_r+\frac{2u}{r})^2 r^2\mathrm{d}r\d\tau\leq E_0,
\end{equation}
and
\begin{equation}\label{3.20}
a(t)\in H^1([0,T]),\ c_0\leq a(t) \leq C_T ,\ t\in(0,T),
\end{equation}
where $E_0:=\int_{\varepsilon}^{a_0}(\frac{1}{2}\rho_0 u_0^2+\frac{1}{\gamma-1}\rho_0^\gamma) r^2\mathrm{d}r,$ $c_0$ and $C_T$ are positive constants. $C_{\gamma}$ is the positive constant defined in Theorem \ref{thm1}
\par
Moreover, if $M<\overline{M}<M_c$, then
\begin{equation}\label{3.18}
\int_{\varepsilon}^{a(t)}(\frac{1}{2}\rho u^2+\f{1}{2(\gamma-1)}\rho^\gamma) r^2\mathrm{d}r+ (\lambda+2\mu)\int_{0}^{t}\int_{\varepsilon}^{a(\tau)}(u_r+\frac{2u}{r})^2 r^2\mathrm{d}r\d\tau\leq E_0.
\end{equation}
\end{lemma}
\pf Multiplying $(\ref{1.2a})_2$ by $u r^2$ and integrating the resulted equation over $(\varepsilon,a(t))$, we obtain after integrating by part and using
(\ref{3.1})-(\ref{3.2}) 
  that
  \begin{equation}\label{2.6}
   \frac{\mathrm{d}}{\mathrm{d} t}\int_{\varepsilon}^{a(t)}(\frac{1}{2}\rho u^2+\frac{1}{\gamma-1}\rho^\gamma) r^2\mathrm{d}r+4\pi\int_{\varepsilon}^{a(t)}\rho u\int_{\varepsilon}^{r}\rho s^2\d s\mathrm{d}r+(\lambda+2\mu)\int_{\varepsilon}^{a(t)}(u_r+\frac{2u}{r})^2 r^2\mathrm{d}r=0.
\end{equation}
The second term on the left hand of \eqref{2.6} can be rewritten as follows
\begin{eqnarray}\label{2.7}
  \nonumber\int_{\varepsilon}^{a(t)}\rho u\int_{\varepsilon}^{r}\rho s^2\d s\mathrm{d}r&=& -\f{1}{2} \int_{\varepsilon}^{a(t)}\rho_t r\int_{\varepsilon}^{r}\rho s^2\d s\mathrm{d}r-\f{1}{2} \int_{\varepsilon}^{a(t)}(\rho u)_r r\int_{\varepsilon}^{r}\rho s^2\d s\mathrm{d}r\\
  \nonumber&=&-\f{1}{2} \f{\d}{\d t}\int_{\varepsilon}^{a(t)}\rho r\int_{\varepsilon}^{r}\rho s^2\d s\mathrm{d}r-\f{1}{2} \int_{\varepsilon}^{a(t)}\rho r\int_{\varepsilon}^{r}(\rho u s^2)_s\d s\mathrm{d}r\\
  \nonumber&&+\f{1}{2}\int_{\varepsilon}^{a(t)}\rho u\int_{\varepsilon}^{r}\rho s^2\d s\mathrm{d}r+\f{1}{2}\int_{\varepsilon}^{a(t)}\rho^2 u r^3\mathrm{d}r\\
  \nonumber&=&\f{1}{2}\int_{\varepsilon}^{a(t)}\rho u\int_{\varepsilon}^{r}\rho s^2\d s\mathrm{d}r-\f{1}{2} \f{\d}{\d t}\int_{\varepsilon}^{a(t)}\rho r\int_{\varepsilon}^{r}\rho s^2\d s\mathrm{d}r.
\end{eqnarray}
which implies
\begin{equation}\label{2.8}
 \int_{\varepsilon}^{a(t)}\rho u\int_{\varepsilon}^{r}\rho s^2\d s\mathrm{d}r = - \f{\d}{\d t}\int_{\varepsilon}^{a(t)}\rho r\int_{\varepsilon}^{r}\rho s^2\d s\mathrm{d}r
\end{equation}
Then, one can deduce from $(\ref{2.6})$ and $(\ref{2.8})$ that
\begin{eqnarray}\label{2.9}
 \nonumber  && \int_{\varepsilon}^{a(t)}(\frac{1}{2}\rho u^2+\frac{1}{\gamma-1}\rho^\gamma) r^2\mathrm{d}r-4\pi\int_{\varepsilon}^{a(t)}\rho r\int_{\varepsilon}^{r}\rho s^2\d s\mathrm{d}r+ (\lambda+2\mu)\int_{0}^{t}\int_{\varepsilon}^{a(\tau)}(u_r+\frac{2u}{r})^2 r^2\mathrm{d}r\d\tau\\
 &\leq& E_0:=\int_{\varepsilon}^{a_0}(\frac{1}{2}\rho_0 u_0^2+\frac{1}{\gamma-1}\rho_0^\gamma) r^2\mathrm{d}r
\end{eqnarray}
We estimate the second term on the left hand side of \eqref{2.9} below. By definition $\Phi_r=-\f{4\pi}{r^2}\int_{\varepsilon}^{r}\rho s^2 \d s$, it holds
\begin{equation}\label{jia2.9}
   4\pi \int_{\varepsilon}^{a(t)}\rho r\int_{\varepsilon}^{r}\rho s^{2}\d s\d r= \f{1}{8\pi}\int_{\varepsilon}^{a(t)}r^2|\Phi_r|^2 \mathrm{d}r+\f{2\pi }{a(t)}(\int_{\varepsilon}^{a(t)}\rho r^{2}\d r)^2.
\end{equation}
Using the elliptic equation $(\ref{1.1})_3$, H\"{o}lder's inequality and interpolation inequality, one can get that
\begin{eqnarray}\label{2.10}
\nonumber && \int_{\varepsilon}^{a(t)}r^2|\Phi_r|^2 \mathrm{d}r=\f{1}{4\pi}\|\nabla\Phi\|^2_{L^2(\Omega^\varepsilon(t))}\leq \f{1}{4\pi}\|\nabla\Phi\|^2_{L^2(\mathbb{R}^3)}\\
&&\leq\|\rho\|_{L^p(\Omega^\varepsilon(t))}\|\Phi\|_{L^{p^\prime}(\mathbb{R}^3)}\leq\|\rho\|^{\theta}_{L^{\gamma}(\Omega^\varepsilon(t))}\|\rho\|^{1-\theta}_{L^{1}(\Omega^\varepsilon(t))}\|\Phi\|_{L^{p^\prime}(\mathbb{R}^3)},
\end{eqnarray}
where $\f{1}{p}+\f{1}{p^\prime}=1,\ \f{1}{p}=\f{\theta}{\gamma}+\f{1-\theta}{1},\ \|\cdot\|_{L^{p}(\Omega^\varepsilon(t))}=(\int_{\Omega^\varepsilon(t)}|\cdot|^p \mathrm{d}\x)^{\f{1}{p}}=(4\pi\int_{\varepsilon}^{a(t)}r^2|\cdot|^p \mathrm{d}r)^{\f{1}{p}}$ and $\Omega^\varepsilon(t):=\{\x\in \mathbb{R}^3|\varepsilon\leq|\x|\leq a(t)\}$ for $t\in~[0,T].$\\
By the elliptic equation $(\ref{1.1})_3$ and the Hardy-Littlewood-Sobolev inequality, it holds that
\begin{equation*}
   \|\Phi\|_{L^{p^\prime}(\mathbb{R}^3)}\leq A_{\gamma}\|\rho\|_{L^{\gamma}(\Omega^\varepsilon(t))},
\end{equation*}
where $\f{1}{p^\prime}=\f{1}{\gamma}-\f{2}{3}$ and $ A_{\gamma}$ is a positive constant just depending on $\gamma$. Thus we have
\begin{equation}\label{2.11}
  \f{1}{8\pi}\int_{\varepsilon}^{a(t)}r^2|\Phi_r|^2 \mathrm{d}r=\f{1}{32\pi^2}\|\nabla\Phi\|^2_{L^2(\Omega^\varepsilon(t))} \leq \f{1}{8\pi} A_{\gamma} \|\rho\|^{1+\theta}_{L^{\gamma}(\Omega^\varepsilon(t))}\|\rho\|^{1-\theta}_{L^{1}(\Omega^\varepsilon(t))},
\end{equation}
where $1+\theta=\f{\gamma}{3(\gamma-1)},\ 1-\theta=\f{5\gamma-6}{3(\gamma-1)}>0$ .\\
Since
 \begin{eqnarray*}
  \f{M}{4\pi}=\int_{\varepsilon}^{a(t)}\rho r^2\mathrm{d}r &\leq& (\int_{\varepsilon}^{a(t)}\rho^\gamma r^2\mathrm{d}r)^{\frac{1}{\gamma}}(\int_{\varepsilon}^{a(t)} r^2\mathrm{d}r)^{1-\frac{1}{\gamma}}
 \leq 3^{\frac{1}{\gamma}-1}(\int_{\varepsilon}^{a(t)}\rho^\gamma r^2\mathrm{d}r)^{\frac{1}{\gamma}}a(t)^{3(1-\frac{1}{\gamma})},
 \end{eqnarray*}
 then we obtain
 \begin{equation}\label{2.120}
   \f{1}{a(t)} \leq [\f{4\pi}{M}3^{\f{1-\gamma}{\gamma}}]^{\f{\gamma}{3(\gamma-1)}}\cdot(\int_{\varepsilon}^{a(t)}\rho^\gamma r^2\mathrm{d}r)^{\frac{1}{3(\gamma-1)}},
 \end{equation}
and
 \begin{equation}\label{2.12}
  \f{2\pi }{a(t)}(\int_{\varepsilon}^{a(t)}\rho r^{2}\d r)^2=\f{2\pi }{a(t)}(\f{M}{4\pi})^2\leq \f{1}{2\sqrt[3]{3}}(4\pi)^{\f{3-2\gamma}{3(\gamma-1)}}M^{\f{5\gamma-6}{3(\gamma-1)}}(\int_{\varepsilon}^{a(t)}\rho^\gamma r^2\mathrm{d}r)^{\f{1}{3(\gamma-1)}}.
 \end{equation}
By \eqref{jia2.9}, \eqref{2.11} and \eqref{2.12}, we have
\begin{eqnarray}\label{2.13}
   \nonumber&&\frac{1}{\gamma-1}\int_{\varepsilon}^{a(t)}\rho^\gamma r^2\mathrm{d}r-4\pi\int_{\varepsilon}^{a(t)}\rho r\int_{\varepsilon}^{r}\rho s^2\d s\mathrm{d}r\\
   \nonumber&=&\frac{1}{\gamma-1}\int_{\varepsilon}^{a(t)}\rho^\gamma r^2\mathrm{d}r-\f{1}{8\pi}\int_{\varepsilon}^{a(t)}r^2|\Phi_r|^2 \mathrm{d}r-\f{2\pi }{a(t)}(\int_{\varepsilon}^{a(t)}\rho r^{2}\d r)^2\\
 &\geq&\frac{1}{\gamma-1}\int_{\varepsilon}^{a(t)}\rho^\gamma r^2\mathrm{d}r-B M^{\f{5\gamma-6}{3(\gamma-1)}}(\int_{\varepsilon}^{a(t)}\rho^\gamma r^2\mathrm{d}r)^{\f{1}{3(\gamma-1)}},
\end{eqnarray}
where $B=(4\pi)^{\f{1}{3(\gamma-1)}}(\f{1}{2\sqrt[3]{3}}(4\pi)^{-\f{2}{3}}+\f{A_{\gamma}}{8\pi}).$\par
We want to show the negative gravitational energy can be dominated by the positive kinetic-internal  energy below that is
\begin{align}\label{f2.110}
 \nonumber&\frac{1}{\gamma-1}\int_{\varepsilon}^{a(t)}\rho^\gamma r^2\mathrm{d}r-4\pi\int_{\varepsilon}^{a(t)}\rho r\int_{\varepsilon}^{r}\rho s^2\d s\mathrm{d}r\\
\nonumber& \geq\frac{1}{\gamma-1}\int_{\varepsilon}^{a(t)}\rho^\gamma r^2\mathrm{d}r-B M^{\f{5\gamma-6}{3(\gamma-1)}}(\int_{\varepsilon}^{a(t)}\rho^\gamma r^2\mathrm{d}r)^{\f{1}{3(\gamma-1)}}\\
& \geq
\left\{ \aligned
  &C_{\gamma}\int_{\varepsilon}^{a(t)}\rho^\gamma r^2\mathrm{d}r,&M<M_c ;\\
  &\f{1}{2(\gamma-1)}\int_{\varepsilon}^{a(t)}\rho^\gamma r^2\mathrm{d}r ,&M<\overline{M}<M_c.
 \endaligned\right.
\end{align}
For the case $\gamma=\f{4}{3}$, by \eqref{2.9} and \eqref{2.13}, we obtain
 \begin{align}\label{f2.11}
  \nonumber E_0&\geq  3\int_{\varepsilon}^{a(t)}\rho^{\f{4}{3}} r^2\mathrm{d}r-BM^{\f{2}{3}}\int_{\varepsilon}^{a(t)}\rho^{\f{4}{3}} r^2\mathrm{d}r\\
   \nonumber &=\f{3}{2}\int_{\varepsilon}^{a(t)}\rho^{\f{4}{3}} r^2\mathrm{d}r+(\f{3}{2}-BM^{\f{2}{3}})\int_{\varepsilon}^{a(t)}\rho^{\f{4}{3}} r^2\mathrm{d}r\\
   &\geq {\left\{ \aligned
   &(3-BM^{\f{2}{3}})\int_{\varepsilon}^{a(t)}\rho^{\f{4}{3}} r^2\mathrm{d}r,~~  M<(\f{3}{B})^{\f{3}{2}};\\
    &\f{3}{2}\int_{\varepsilon}^{a(t)}\rho^{\f{4}{3}} r^2\mathrm{d}r, ~~M\leq(\f{3}{2B})^{\f{3}{2}}<(\f{3}{B})^{\f{3}{2}}.
 \endaligned\right.
 }
 \end{align}
For the case $\f{6}{5}<\gamma<\f{4}{3},$  define
\begin{equation}\label{f2.1}
  f(s)=\frac{1}{\gamma-1}s-BM^{\f{5\gamma-6}{3(\gamma-1)}}s^{\f{1}{3(\gamma-1)}},~s>0.
\end{equation}
Then, we have
\begin{align}\label{f2.2}
   &f^\prime(s)=\frac{1}{\gamma-1}-\f{1}{3(\gamma-1)}BM^{\f{5\gamma-6}{3(\gamma-1)}}s^{\f{4-3\gamma}{3(\gamma-1)}},\\
   &f^{\prime\prime}(s)=-\f{1}{3(\gamma-1)}\f{4-3\gamma}{3(\gamma-1)}BM^{\f{5\gamma-6}{3(\gamma-1)}}s^{\f{7-6\gamma}{3(\gamma-1)}},\\
   &f^{\prime\prime}(s)<0,~s>0.
\end{align}
By \eqref{2.9}, \eqref{2.13} and \eqref{f2.1}, it holds that
\begin{equation}\label{fp}
   f(\int_{\varepsilon}^{a(t)}\rho^\gamma r^2\mathrm{d}r)\leq E_0.
\end{equation}
Define
\begin{equation}\label{f2.2s}
s^*=(\f{B}{3})^{-\f{3(\gamma-1)}{4-3\gamma}}M^{-\f{5\gamma-6}{4-3\gamma}},
\end{equation}
which solves
$ f^\prime(s^*)=0,$
and then there is
\begin{equation}\label{f2.3}
  f(s^*)=\f{4-3\gamma}{\gamma-1}(\f{B}{3})^{-\f{3(\gamma-1)}{4-3\gamma}}M^{-\f{5\gamma-6}{4-3\gamma}}.
\end{equation}
Thus, $f(s)$ increases strictly in $(0,s^*)$ and decreases strictly in $(s^*,+\infty).$\\
We claim that under the condition
\begin{equation}\label{f2.41}
    M<M_c=[\f{4-3\gamma}{\gamma-1}(\f{B}{3})^{-\f{3(\gamma-1)}{4-3\gamma}}]^{\f{4-3\gamma}{5\gamma-6}}E_0^{-\f{4-3\gamma}{5\gamma-6}},
\end{equation}
  it holds that $\int_{\varepsilon}^{a(t)}\rho^\gamma r^2\mathrm{d}r<s^*$ and then $f(\int_{\varepsilon}^{a(t)}\rho^\gamma r^2\mathrm{d}r)>\frac{4-3\gamma}{\gamma-1}\int_{\varepsilon}^{a(t)}\rho^\gamma r^2\mathrm{d}r$.
Indeed, if $M<M_c$, then by direct calculation, it holds that $f(s^*)>E_0$ and
 \begin{align*}
   s^*>(\f{B}{3})^{-\f{3(\gamma-1)}{4-3\gamma}}\f{\gamma-1}{4-3\gamma}(\f{B}{3})^{\f{3(\gamma-1)}{4-3\gamma}}E_0>\f{\gamma-1}{4-3\gamma}E_0>2(\gamma-1)E_0,
   ~~\f{6}{5}<\gamma<\f{4}{3}.
\end{align*}
Since
\begin{equation}\label{f2.60}
    \int_{\varepsilon}^{a_0}\rho_0^\gamma r^2\mathrm{d}r\leq (\gamma-1)E_0< s^*,
\end{equation}
\begin{equation}\label{f2.600}
    f(\int_{\varepsilon}^{a(t)}\rho^\gamma r^2\mathrm{d}r)\leq E_0< f(s^*),
\end{equation}
and the continuity of $\int_{\varepsilon}^{a(t)}\rho^\gamma r^2\mathrm{d}r$ with respect to $t$, we conclude that
\begin{align}\label{f2.6}
    \int_{\varepsilon}^{a(t)}\rho^\gamma r^2\mathrm{d}r<s^*.
\end{align}
If \eqref{f2.6} is not correct, then by the continuity of $\int_{\varepsilon}^{a(t)}\rho^\gamma r^2\mathrm{d}r$ with respect to $t$ and \eqref{f2.60}, there exists $t_0>0$ such that
$\int_{\varepsilon}^{a(t_0)}\rho^\gamma r^2\mathrm{d}r=s^*$, and then
\begin{equation*}
  f(\int_{\varepsilon}^{a(t_0)}\rho^\gamma r^2\mathrm{d}r)=f(s^*)> E_0,
\end{equation*}
which contradicts with \eqref{f2.600}. Therefore, \eqref{f2.6} holds for $t\geq0$.
By \eqref{2.13}, \eqref{fp} and \eqref{f2.6}, we obtain
\begin{align}\label{f2.91}
 \nonumber E_0&\geq \frac{1}{\gamma-1}\int_{\varepsilon}^{a(t)}\rho^\gamma r^2\mathrm{d}r-B M^{\f{5\gamma-6}{3(\gamma-1)}} (\int_{\varepsilon}^{a(t)}\rho^\gamma r^2\mathrm{d}r)^{\f{1}{3(\gamma-1)}}\\
  \nonumber&\geq \frac{1}{\gamma-1}\int_{\varepsilon}^{a(t)}\rho^\gamma r^2\mathrm{d}r- B M^{\f{5\gamma-6}{3(\gamma-1)}} (s^*)^{\f{4-3\gamma}{3(\gamma-1)}} \int_{\varepsilon}^{a(t)}\rho^\gamma r^2\mathrm{d}r \\
 \nonumber&=\frac{4-3\gamma}{\gamma-1}\int_{\varepsilon}^{a(t)}\rho^\gamma r^2\mathrm{d}r.
\end{align}
We claim that under the condition
\begin{equation}\label{f2.4}
    M<\overline{M}:=[\f{4-3\gamma}{\gamma-1}(\f{B}{3})^{-\f{3(\gamma-1)}{4-3\gamma}}]^{\f{4-3\gamma}{5\gamma-6}}(lE_0)^{-\f{4-3\gamma}{5\gamma-6}}<M_c,
\end{equation}
for some $l> 1$, it holds that $f(\int_{\varepsilon}^{a(t)}\rho^\gamma r^2\mathrm{d}r)>\frac{1}{2(\gamma-1)}\int_{\varepsilon}^{a(t)}\rho^\gamma r^2\mathrm{d}r$.
Indeed, if $M<\overline{M}<M_c$, then it holds that $f(s^*)>lE_0$ and
\begin{align*}
   s^*&= (\f{B}{3})^{-\f{3(\gamma-1)}{4-3\gamma}}M^{-\f{5\gamma-6}{4-3\gamma}}>\f{\gamma-1}{4-3\gamma}lE_0>2l(\gamma-1)E_0.
\end{align*}
Choose a positive constant $0<\alpha<1$, such that $\alpha l\geq 1$ and then
\begin{align}\label{f2.7}
  \nonumber f(\alpha s^*)&=\frac{1}{\gamma-1}\alpha s^*-BM^{\f{5\gamma-6}{3(\gamma-1)}}(\alpha s^*)^{\f{1}{3(\gamma-1)}}\\
  \nonumber&=\alpha(\frac{1}{\gamma-1}s^*-BM^{\f{5\gamma-6}{3(\gamma-1)}}(s^*)^{\f{1}{3(\gamma-1)}})
  +\alpha(1-\alpha^{\f{4-3\gamma}{3(\gamma-1)}})BM^{\f{5\gamma-6}{3(\gamma-1)}}(s^*)^{\f{1}{3(\gamma-1)}}\\
  &\geq \alpha f( s^*)>\alpha lE_0 \geq E_0\geq f(\int_{\varepsilon}^{a(t)}\rho^\gamma r^2\mathrm{d}r).
\end{align}
Based on \eqref{f2.6}, \eqref{f2.7} and the growth of $f(s)$ in $(0,s^*)$, we obtain
\begin{equation}\label{f2.8}
    \int_{\varepsilon}^{a(t)}\rho^\gamma r^2\mathrm{d}r<\alpha s^*.
\end{equation}
By \eqref{2.13}, \eqref{f2.7} and \eqref{f2.8}, we obtain
\begin{align}\label{f2.9}
  \nonumber E_0&\geq \frac{1}{\gamma-1}\int_{\varepsilon}^{a(t)}\rho^\gamma r^2\mathrm{d}r-B M^{\f{5\gamma-6}{3(\gamma-1)}} (\int_{\varepsilon}^{a(t)}\rho^\gamma r^2\mathrm{d}r)^{\f{1}{3(\gamma-1)}}\\
  \nonumber&\geq \frac{1}{\gamma-1}\int_{\varepsilon}^{a(t)}\rho^\gamma r^2\mathrm{d}r- B M^{\f{5\gamma-6}{3(\gamma-1)}}(\alpha s^*)^{\f{4-3\gamma}{3(\gamma-1)}} \int_{\varepsilon}^{a(t)}\rho^\gamma r^2\mathrm{d}r \\
&=\frac{1}{\gamma-1}\int_{\varepsilon}^{a(t)}\rho^\gamma r^2\mathrm{d}r- 3\alpha ^{\f{4-3\gamma}{3(\gamma-1)}} \int_{\varepsilon}^{a(t)}\rho^\gamma r^2\mathrm{d}r.
\end{align}
If choose $l,\alpha$ satisfying $0<\alpha<1$, $~l>1,~\alpha l\geq1,~\frac{1}{2}- \alpha ^{\f{4-3\gamma}{3(\gamma-1)}}\geq 0$, then
\begin{align}\label{}
  \nonumber&\frac{1}{\gamma-1}\int_{\varepsilon}^{a(t)}\rho^\gamma r^2\mathrm{d}r- 3\alpha ^{\f{4-3\gamma}{3(\gamma-1)}} \int_{\varepsilon}^{a(t)}\rho^\gamma r^2\mathrm{d}r \\ \nonumber&=\frac{1}{2(\gamma-1)}\int_{\varepsilon}^{a(t)}\rho^\gamma r^2\mathrm{d}r+(\frac{1}{2(\gamma-1)}- 3\alpha ^{\f{4-3\gamma}{3(\gamma-1)}} ) \int_{\varepsilon}^{a(t)}\rho^\gamma r^2\mathrm{d}r\\
 \nonumber &\geq \frac{1}{2(\gamma-1)}\int_{\varepsilon}^{a(t)}\rho^\gamma r^2\mathrm{d}r+3(\frac{1}{2}- \alpha ^{\f{4-3\gamma}{3(\gamma-1)}} ) \int_{\varepsilon}^{a(t)}\rho^\gamma r^2\mathrm{d}r\\
  &\geq\frac{1}{2(\gamma-1)}\int_{\varepsilon}^{a(t)}\rho^\gamma r^2\mathrm{d}r.
\end{align}
 Thus, one obtains \eqref{3.19} and \eqref{3.18} from \eqref{2.9} and \eqref{f2.110}.
 By \eqref{3.19}, \eqref{2.120} and the fact that
 \begin{eqnarray}\label{a0}
  \nonumber\int_{0}^{t} \int_{\varepsilon}^{a(\tau)}(u_{r}+\f{2u}{r})^2r^2\d r\d \tau &=&\int_{0}^{t} \int_{\varepsilon}^{a(\tau)}(u_{r}^{2}r^2+2u^2+2(u^2 r)_r)\d r\d \tau \\
  &=& \int_{0}^{t} \int_{\varepsilon}^{a(\tau)}(u_{r}^{2}+\f{2u^2}{r^2})r^2\d r\d \tau+2\int_{0}^{t}a(\tau)(a^\prime(\tau))^2 \d \tau,
\end{eqnarray}
we obtain \eqref{3.20}. The proof of Lemma \ref{lem2.1} is completed.\\

In order to use the similar method to prove the global existence as in \cite{JZ2001,SG}, we need the uniform a-priori estimate around the symmetry center below.
\begin{lemma}\label{lem4.1}
  Let $T>0$, $\f{6}{5}<\gamma\leq\f{4}{3}$ and $(\rho,u,a)$  be any strong solution to the FBVP \eqref{1.2a}-\eqref{3.2} for $t\in[0,T]$ under the same assumption mass $M<M_c$, then there exists a positive constant $ C_{0,T}>0$ depending on $E_0,~M$ and $T$, but independent of $\varepsilon$, such that
  \begin{equation}\label{4.1}
    \int_{0}^{T}\int_{\varepsilon}^{a(t)}\rho^{2\gamma}(r,t)r^{12}\mathrm{d}r\mathrm{d}t\leq C_{0,T}.
  \end{equation}
\end{lemma}
\proof
We multiply $(\ref{1.2a})_2$ by $\varphi(r)(=r^3)$ and integrate the resulted equation over $(r,a(t))\ (r\in[\varepsilon ,a(t)])$ to obtain
\begin{align}
  \nonumber&(\rho u^2+\rho ^\gamma-(\lambda+2\mu) u_r-(\lambda+2\mu) \frac{2u}{r}) \varphi \\
\nonumber=& \partial_t \int_{r}^{a(t)}  \rho u \varphi \mathrm{d}y-\int_{r}^{a(t)}(\rho u^2+\rho ^\gamma-(\lambda+2\mu) u_y
   -(\lambda+2\mu) \frac{2u}{y})\varphi_y\mathrm{d}y\\
 & +\int_{r}^{a(t)}\frac{2\rho u^2}{y}\varphi \mathrm{d}y+\int_{r}^{a(t)}\frac{4\pi\rho \varphi }{y^2} \int_{\varepsilon}^{y}\rho s^2\d s\mathrm{d}y,\label{4.2}
\end{align}
which yields
 \begin{eqnarray}\label{4.3}
     \nonumber\rho ^{2\gamma}\varphi &=& (\lambda+2\mu)\rho ^\gamma(u_r+\frac{2u}{r})\varphi-\rho ^{1+\gamma}u^2\varphi+\rho ^\gamma \partial_t \int_{r}^{a(t)}  \rho u \varphi \mathrm{d}y-\rho ^\gamma \int_{r}^{a(t)}(\rho u^2+\rho ^\gamma)\varphi_y\mathrm{d}y \\
     & & +(\lambda+2\mu)\rho ^\gamma \int_{r}^{a(t)}( u_y+\frac{2u}{y})\varphi_y\mathrm{d}y+\rho ^\gamma \int_{r}^{a(t)}\frac{2\rho u^2}{y}\varphi \mathrm{d}y+\rho^\gamma\int_{r}^{a(t)}\frac{4\pi\rho \varphi }{y^2} \int_{\varepsilon}^{y}\rho s^2\d s\mathrm{d}y.
  \end{eqnarray}
It follows from $(\ref{1.2a})_1$ that $\rho ^\gamma$ satisfies
  \begin{equation}\label{4.4}
    (\rho ^\gamma)_t+(\rho ^\gamma u)_r+\frac{2\gamma \rho ^\gamma u}{r}=(1-\gamma)\rho ^\gamma u_r.
  \end{equation}
Then, one has
 \begin{eqnarray}\label{4.5}
 \nonumber &&\rho ^\gamma \partial_t \int_{r}^{a(t)}  \rho u \varphi \mathrm{d}y=\partial_t(\rho ^\gamma \int_{r}^{a(t)}  \rho u \varphi \mathrm{d}y)-\partial_t\rho ^\gamma\int_{r}^{a(t)}  \rho u \varphi \mathrm{d}y=\partial_t(\rho ^\gamma \int_{r}^{a(t)}  \rho u \varphi \mathrm{d}y) \\
  &&+\partial_r(\rho ^\gamma u \int_{r}^{a(t)}\varphi\rho u \mathrm{d}y) +\{\frac{2\gamma \rho ^\gamma u}{r}+(\gamma-1)\rho ^\gamma u_r\}\int_{r}^{a(t)}  \rho u \varphi \mathrm{d}y+\varphi\rho ^{1+\gamma}u^2.
 \end{eqnarray}
Substituting (\ref{4.5}) into (\ref{4.3}), we get
 \begin{eqnarray}\label{4.6}
 \nonumber \rho ^{2\gamma}\varphi &=& (\lambda+2\mu)\rho ^\gamma(u_r+\frac{2u}{r})\varphi+\partial_t(\rho ^\gamma \int_{r}^{a(t)}  \rho u \varphi \mathrm{d}y)+\partial_r(\rho ^\gamma u \int_{r}^{a(t)}\rho u\varphi \mathrm{d}y)\ \ \ \ \ \  \\
     \nonumber& &+\{\frac{2\gamma \rho ^\gamma u}{r}+(\gamma-1)\rho ^\gamma u_r\}\int_{r}^{a(t)}  \rho u \varphi \mathrm{d}y-\rho ^\gamma \int_{r}^{a(t)}(\rho u^2+\rho ^\gamma)\varphi_y\mathrm{d}y\\
     & &+(\lambda+2\mu)\rho ^\gamma \int_{r}^{a(t)}( u_y+\frac{2u}{y})\varphi_y\mathrm{d}y+\rho ^\gamma \int_{r}^{a(t)}\frac{2\rho u^2}{y}\varphi \mathrm{d}y +\rho^\gamma\int_{r}^{a(t)}\frac{4\pi\rho \varphi }{y^2} \int_{\varepsilon}^{y}\rho s^2\d s\mathrm{d}y.
  \end{eqnarray}
 Multiplying (\ref{4.6}) by $\varphi ^3$ and integrating over $[\varepsilon,a(t)]\times[0,T]$ lead to
\begin{equation}\label{4.7}
    \int_{0}^{T}\int_{\varepsilon}^{a(t)}\rho ^{2\gamma}\varphi^4 \mathrm{d}r\mathrm{d}t  = \int_{0}^{T}\int_{\varepsilon}^{a(t)}\{R.H.S\  of \ (\ref{4.6})\}\varphi ^3 \mathrm{d}r\mathrm{d}t=\sum_{i=1}^{6} I_i .
\end{equation}
The right hand side terms of \eqref{4.7} can be estimated as follows:
\begin{eqnarray*}
  |I_1| &=& |(\lambda+2\mu)\int_{0}^{T}\int_{\varepsilon}^{a(t)}\rho ^\gamma (u_r+\frac{2u}{r})\varphi^4 \mathrm{d}r\mathrm{d}t|\\
  &\leq& \delta\int_{0}^{T}\int_{\varepsilon}^{a(t)}\rho ^{2\gamma} \varphi^4 \mathrm{d}r\mathrm{d}t+C\delta^{-1}\int_{0}^{T}\int_{\varepsilon}^{a(t)} (u_r
  ^2+\frac{4u^2}{r^2})\varphi^4 \mathrm{d}r\mathrm{d}t \\
  &\leq&  \delta\int_{0}^{T}\int_{\varepsilon}^{a(t)}\rho ^{2\gamma} \varphi^4 \mathrm{d}r\mathrm{d}t+\delta^{-1}C_{0,T} .
\end{eqnarray*}
It is easy to see that by (\ref{3.19}) and (\ref{3.20}):
\begin{equation}\label{4.8}
    |\int_{r}^{a(t)}  \rho u \varphi \mathrm{d}y|\leq\frac{1}{2}\int_{r}^{a(t)} ( \rho u^2+\rho )\varphi \mathrm{d}y \leq C_{0,T},
\end{equation}
hence,
\begin{eqnarray*}
  |I_2| &=& |\int_{0}^{T}\int_{\varepsilon}^{a(t)}\varphi^3\partial_t(\rho ^\gamma \int_{r}^{a(t)}  \rho u \varphi \mathrm{d}y) \mathrm{d}r\mathrm{d}t|
\leq C\sup_{0\leq t\leq T}|\int_{\varepsilon}^{a(t)}\rho ^\gamma\varphi^3( \int_{r}^{a(t)}  \rho u \varphi \mathrm{d}y) \mathrm{d}r\mathrm{d}t|\\
     &\leq& C_{0,T} \sup_{0\leq t\leq T}\int_{\varepsilon}^{a(t)}\rho^\gamma \varphi^3\mathrm{d}r \leq C_{0,T},
\end{eqnarray*}
\begin{eqnarray*}
     |I_3| &=& |\int_{0}^{T}\int_{\varepsilon}^{a(t)}\partial_r(\rho ^\gamma u \int_{r}^{a(t)}\varphi\rho u \mathrm{d}y)\varphi^3 \mathrm{d}r\mathrm{d}t| =|\int_{0}^{T}\int_{\varepsilon}^{a(t)}3\rho ^\gamma u\varphi^2\varphi_r ( \int_{r}^{a(t)}\varphi\rho u \mathrm{d}y) \mathrm{d}r\mathrm{d}t| \\
     &\leq& \delta\int_{0}^{T}\int_{\varepsilon}^{a(t)}\rho^{2\gamma}\varphi^4\mathrm{d}r\mathrm{d}t
     +\delta^{-1}C_{0,T}\int_{0}^{T}\int_{\varepsilon}^{a(t)}u^2 \varphi_r^2\mathrm{d}r\mathrm{d}t \\
     &\leq& \delta\int_{0}^{T}\int_{\varepsilon}^{a(t)}\rho^{2\gamma}\varphi^4\mathrm{d}r\mathrm{d}t+ \delta^{-1}C_{0,T},
 \end{eqnarray*}
 \begin{eqnarray*}
     |I_4| &=&|\int_{0}^{T}\int_{\varepsilon}^{a(t)}\varphi^3\{\frac{2\gamma \rho ^\gamma u}{r}+(\gamma-1)\rho ^\gamma u_r\}(\int_{r}^{a(t)}  \rho u \varphi \mathrm{d}y )\mathrm{d}r\mathrm{d}t |  \\
     &\leq& C_{0,T} \int_{0}^{T}\int_{\varepsilon}^{a(t)}\varphi^3 |\frac{2\gamma \rho ^\gamma u}{r}+(\gamma-1)\rho ^\gamma u_r|\mathrm{d}r\mathrm{d}t \\
      &\leq& \delta\int_{0}^{T}\int_{\varepsilon}^{a(t)}\rho^{2\gamma}\varphi^4\mathrm{d}r\mathrm{d}t
      +\delta^{-1}C_{0,T} \int_{0}^{T}\int_{\varepsilon}^{a(t)}(\frac{u^2}{r^2}+u_r^2)\varphi^2\mathrm{d}r\mathrm{d}t \\
     &\leq& \delta\int_{0}^{T}\int_{\varepsilon}^{a(t)}\rho^{2\gamma}\varphi^4\mathrm{d}r\mathrm{d}t
      +\delta^{-1}C_{0,T}.
  \end{eqnarray*}
 By (\ref{3.19}) and (\ref{3.20}), we get
  \begin{equation}\label{4.9}
    |\int_{r}^{a(t)}(\rho u^2+\rho ^\gamma))\varphi_y\mathrm{d}y|\leq CE_{0},\ |\int_{r}^{a(t)}\frac{2\rho u^2}{y}\varphi \mathrm{d}y|\leq CE_{0},
  \end{equation}
  \begin{equation}\label{4.10}
    |\int_{0}^{T}\int_{r}^{a(t)}( u_y+\frac{2u}{y})\varphi_y\mathrm{d}y\d t|\leq C\int_{0}^{T} a(t)\d t+C\int_{0}^{T}\int_{\varepsilon}^{a(t)}( u_y^2+\frac{4u^2}{y^2})\varphi_y^2\mathrm{d}y\d t\leq C_{0,T},
  \end{equation}
 which gives
\begin{eqnarray*}
\nonumber |I_5| &=& |\int_{0}^{T}\int_{\varepsilon}^{a(t)}\rho ^\gamma\varphi^3\{\int_{r}^{a(t)}(\rho u^2+\rho ^\gamma))\varphi_y\mathrm{d}y\\
\nonumber&&+(\lambda+2\mu)\int_{r}^{a(t)}( u_y+\frac{2u}{y})\varphi_y\mathrm{d}y+\int_{r}^{a(t)}\frac{2\rho u^2}{y}\varphi \mathrm{d}y\}\mathrm{d}r\mathrm{d}t| \\
    &\leq& C_{0,T}(\int_{0}^{T}\int_{\varepsilon}^{a(t)}\rho^\gamma\varphi^3\mathrm{d}r\mathrm{d}t
    +\sup_{0\leq t\leq T}\int_{\varepsilon}^{a(t)}\rho^\gamma \varphi^3\mathrm{d}r)
    \leq C_{0,T}.
\end{eqnarray*}
\begin{align*}
   \nonumber |I_6| &= |\int_{0}^{T}\int_{\varepsilon}^{a(t)}\rho^\gamma\varphi^3\int_{r}^{a(t)}\frac{4\pi\rho \varphi }{y^2} \int_{\varepsilon}^{y}\rho s^2\d s\mathrm{d}y\d r\d t|\\
   &\leq CM^2\int_{0}^{T}\int_{\varepsilon}^{a(t)}\rho^\gamma\varphi^3r^{-1}\d r\d t\leq C_{0,T}.
\end{align*}
 We finally get
 \begin{equation*}
    \int_{0}^{T}\int_{\varepsilon}^{a(t)}\rho ^{2\gamma}\varphi^4\mathrm{d}r\mathrm{d}t\leq 3\delta\int_{0}^{T}\int_{\varepsilon}^{a(t)}\rho ^{2\gamma}\varphi^4\mathrm{d}r\mathrm{d}t+\delta^{-1}C_{0,T},
 \end{equation*}
 Choosing $\delta=\frac{1}{6}$ , we obtain
 \begin{equation*}
    \int_{0}^{T}\int_{\varepsilon}^{a(t)}\rho ^{2\gamma}\varphi^4\mathrm{d}r\mathrm{d}t\leq C_{0,T}.
 \end{equation*}

Use the same method as the proof of Lemma 3.4 in \cite{GLX2012}, we gets the bounds of particle path which will be used to make other estimates below.
\begin{lemma}\label{lem2.2}
  Let $T>0$, $\f{6}{5}<\gamma\leq\f{4}{3}$ and $(\rho,u,a)$  be any strong solution to the FBVP \eqref{1.2a}-\eqref{3.2} for $t\in[0,T]$ under the same assumption mass $M<M_c$. Then
  \begin{equation}\label{3.23}
 (E_0/C_\gamma)^{-\frac{1}{3(\gamma-1)}}x^{\frac{\gamma}{3(\gamma-1)}} \leq r(x,\tau)\leq a(\tau),\ (x,\tau)\in[0,\f{M}{4\pi}]\times[0,T],
  \end{equation}
  \begin{equation}\label{3.24}
     (E_0/C_\gamma)^{-\frac{1}{\gamma-1}}(x_2-x_1)^{\frac{\gamma}{\gamma-1}} \leq r^3(x_2,\tau)- r^3(x_1,\tau),\ 0\leq x_1  < x_2\leq \f{M}{4\pi},\ \tau\in[0,T].
  \end{equation}
\end{lemma}

\begin{lemma}\label{lem2.3}
  Let $T>0$, $\f{6}{5}<\gamma\leq\f{4}{3}$ and $(\rho,u,a)$  be any strong solution to the FBVP \eqref{1.2a}-\eqref{3.2} for $\tau\in[0,T]$. If the mass $M<M_c$, then
\begin{equation}\label{3.25}
\rho_0(r(0))e^{-c_{x,T}/(\lambda+2\mu)}\leq \rho (r(t),t)\leq \rho_0(r(0))e^{C_{x,T}/(\lambda+2\mu)},\  \forall t\in [0,T],
\end{equation}
where $r(t)$ is particle path defined as \eqref{3.30} and $C_{x,T},~c_{x,T}$ are positive constants.
\end{lemma}
\proof Define
 $$
 \xi=\int_{a(t)}^{r}\rho u \mathrm{d}y,\quad \eta=\rho u^2(r,t)-\rho u^2(a(t),t)+\int_{a(t)}^{r}\frac{2\rho u^2}{y}\mathrm{d}y.
 $$
A direct calculation together with (\ref{1.2a}) and \eqref{3.2} gives rise to
\begin{eqnarray}\label{3.26}
  \nonumber\xi_t+\eta+F &=& \int_{a(t)}^{r}(\rho u )_t\mathrm{d}y +\rho u^2(r,t)-2\rho u^2(a(t),t)+\int_{a(t)}^{r}\frac{2\rho u^2}{y}\mathrm{d}y\\
 \nonumber  &&+(\rho^\gamma-(\lambda+2\mu) u_r-(\lambda+2\mu)\frac{2u}{r})(r,t)\\
  \nonumber &=& \int_{r}^{a(t)}(\rho u^2+\rho^\gamma-(\lambda+2\mu) u_y-(\lambda+2\mu)\frac{2u}{y})_y\mathrm{d}y+\int_{r}^{a(t)}\f{4\pi\rho}{y^2}\int_{\varepsilon}^{r}\rho s^2\d s\d y\\
  \nonumber &&+(\rho u^2+\rho^\gamma-(\lambda+2\mu) u_r-(\lambda+2\mu)\frac{2u}{r})(r,t)-2\rho u^2(a(t),t)\\
   &=& -\rho u^2(a(t),t)+\int_{r}^{a(t)}\f{4\pi\rho}{y^2}\int_{\varepsilon}^{y}\rho s^2\d s\d y.
\end{eqnarray}
Rewrite $(\ref{1.2a})_1$ as
$$\rho_t+\rho_r u+\rho (u_r+\frac{2u}{r})=0,$$
which together with $(\ref{1.2a})_2$, yields
\begin{equation}\label{3.27}
  ((\lambda+2\mu)\ln\rho)_t+((\lambda+2\mu)\ln\rho)_r u+\rho^\gamma-F=0.
\end{equation}
It follows from (\ref{3.26}) and (\ref{3.27}) that
\begin{eqnarray*}
  &&(\xi+(\lambda+2\mu)\ln\rho)_t+(\xi+(\lambda+2\mu)\ln\rho)_r u+\rho^\gamma\\
  &=& -\rho u^2(a(t),t)-\eta+u \xi_r+\int_{r}^{a(t)}\f{4\pi\rho}{y^2}\int_{\varepsilon}^{y}\rho s^2\d s\d y \\
   &=& \int_{r}^{a(t)} \frac{2\rho u^2}{y}\mathrm{d}y+\int_{r}^{a(t)}\f{4\pi\rho}{y^2}\int_{\varepsilon}^{y}\rho s^2\d s\d y,
\end{eqnarray*}
Thus, one has
\begin{equation}\label{3.28}
    \frac{\mathrm{D}}{\mathrm{D}t}(\xi+(\lambda+2\mu)\ln\rho)+\rho^\gamma=\int_{r}^{ a(t)} \frac{2\rho u^2}{y}\mathrm{d}y+\int_{r}^{a(t)}\f{4\pi\rho}{y^2}\int_{\varepsilon}^{y}\rho s^2\d s\d y,
\end{equation}
where $\frac{\mathrm{D}}{\mathrm{D}t}=\partial_t+u \partial_r.$
Integrating (\ref{3.28}) with respect to time $t$ shows
\begin{align}
   & \nonumber(\xi+(\lambda+2\mu)\ln\rho)(r(t),t)+\int_{0}^{t}\rho^\gamma(r(s),s)\mathrm{d}s
\\ \label{3.29}
  =&(\xi+(\lambda+2\mu)\ln\rho)(r(0),0) +\int_{0}^{t}\int_{r(\tau)}^{a(\tau)}\frac{2\rho u^2}{y}\mathrm{d}y\mathrm{d}\tau+\int_{0}^{t}\int_{r(\tau)}^{a(\tau)}\f{4\pi\rho}{y^2}\int_{\varepsilon}^{y}\rho s^2\d s\d y\d\tau,
\end{align}
where $r(t)$ is particle path defined as
\begin{eqnarray}\label{3.30}
  \left\{%
  \begin{array}{l}
  \frac{\mathrm{d}}{\mathrm{d}t}r(t)=u(r(t),t),
  \\
  r(0)=r,\ r\in[\varepsilon,a_0].
\end{array}%
\right.
\end{eqnarray}
Therefore, it holds that
\begin{eqnarray}\label{3.31}
  \nonumber&&(\lambda+2\mu)\ln\f{\rho(r(t),t)}{\rho_0(r(0))}+\int_{0}^{t}\rho^\gamma(r(s),s)\mathrm{d}s\\
  &&=\int_{r(t)}^{a(t)}\rho u\mathrm{d}y+\int_{0}^{t}\int_{r(\tau)}^{a(\tau)}\frac{2\rho u^2}{y}\mathrm{d}y\mathrm{d}s-\int_{r(0)}^{a_0}\rho_0 u_0\mathrm{d}y+\int_{0}^{t}\int_{r(\tau)}^{a(\tau)}\f{4\pi\rho}{y^2}\int_{\varepsilon}^{y}\rho s^2\d s\d y\d\tau.
\end{eqnarray}
It follows from H\"{o}lder's inequality, \eqref{3.19}, \eqref{3.20} and \eqref{3.23} that
\begin{equation}\label{3.32}
   \int_{r(t)}^{a(t)}\rho u\mathrm{d}y \leq \frac{1}{r^2(x,t)}(\int_{\varepsilon}^{a(t)}\rho u^2 y^2\mathrm{d}y )^{\frac{1}{2}}(\int_{\varepsilon}^{a(t)}\rho y^2\mathrm{d}y )^{\frac{1}{2}}
  \le  (E_0/C_\gamma)^{\frac{2}{3(\gamma-1)}}x^{-\frac{2\gamma}{3(\gamma-1)}} (\f{M E_0}{\pi})^{\frac12},
\end{equation}
\begin{equation}\label{3.33}
    \int_{0}^{t}\int_{r(\tau)}^{a(\tau)}\frac{2\rho u^2}{y}\mathrm{d}y\mathrm{d}\tau \leq \int_{0}^{T}\frac{1}{r^3(x,\tau)}\int_{\varepsilon}^{a(\tau)}2\rho u^2 y^2\mathrm{d}y\mathrm{d}\tau
   \leq  4T(E_0/C_\gamma)^{\frac{1}{\gamma-1}}x^{-\frac{\gamma}{\gamma-1}} E_0.
\end{equation}
and,
\begin{equation}\label{a3.33}
\int_{0}^{t}\int_{r(\tau)}^{a(\tau)}\f{4\pi\rho}{y^2}\int_{\varepsilon}^{y}\rho s^2\d s\d y\d\tau\leq\int_{0}^{t}\f{ 4\pi}{r^4(x,\tau)}\int_{r(\tau)}^{a(\tau)}\rho y^2\int_{\varepsilon}^{y}\rho s^2\d s\d y\d\tau\leq \f{M^2}{4\pi}T (E_0/C_\gamma)^{\frac{4}{3(\gamma-1)}}x^{-\frac{4\gamma}{3(\gamma-1)}}.
\end{equation}
It holds that for the initial data
\begin{equation}\label{3.34}
\int_{r}^{a_0}\rho_0 u_0\mathrm{d}y \leq \frac{1}{r^2}(\int_{\varepsilon}^{a_0}\rho_0 u_0^2 y^2\mathrm{d}y )^{\frac{1}{2}}(\int_{\varepsilon}^{a_0}\rho_0 y^2\mathrm{d}y )^{\frac{1}{2}}\leq (E_0/C_\gamma)^{\frac{2}{3(\gamma-1)}}x^{-\frac{2\gamma}{3(\gamma-1)}} (\f{M E_0}{\pi})^{\frac12}.
\end{equation}
From the above estimates, we have
\begin{equation*}
\rho_0(r(0))e^{-\f{c_{x,T}}{\lambda+2\mu}}\leq \rho (r(t),t)\leq \rho_0(r(0))e^{\f{C_{x,T}}{\lambda+2\mu}},\  \forall t\in [0,T],
\end{equation*}
where $C_{x,T}:=2(C_\gamma)^{-\frac{2}{3(\gamma-1)}} (M /\pi)^{\frac12} E_0^{\frac{3\gamma+1}{6(\gamma-1)}}x^{-\frac{2\gamma}{3(\gamma-1)}} +4T(E_0/C_\gamma)^{\frac{\gamma}{\gamma-1}}x^{-\frac{\gamma}{\gamma-1}}+\f{M^2}{4\pi}T (E_0/C_\gamma)^{\frac{4}{3(\gamma-1)}}x^{-\frac{4\gamma}{3(\gamma-1)}}$ and $c_{x,T}:=C_{x,T}+T\|\rho_0\|^{\gamma}_{L^\infty}e^{\gamma\f{C_{x,T}}{\lambda+2\mu}}.$\\
\subsection{Estimates in Lagrangian coordinates}
The basic energy estimate can be written in Lagrangian coordinates as below.
 \begin{lemma}\label{lem2.10}
  Let $T>0$, $\f{6}{5}<\gamma\leq\f{4}{3}$ and $(\rho,u,a)$  be any strong solution to the FBVP \eqref{1.2a}-\eqref{3.2} for $\tau\in[0,T]$ under the same assumption mass $M<M_c$. Then
  \begin{equation}\label{3.19a}
 \int_{0}^{\f{M}{4\pi}}(\frac{1}{2}u^2+C_\gamma\rho^{\gamma-1})\mathrm{d}x+(\lambda+2\mu)\int_{0}^{\tau}\int_{0}^{\f{M}{4\pi}}
    \{\rho r^4 |u_x|^2+\frac{2u^2}{\rho r^2}\}\mathrm{d}x\mathrm{d}s+(\lambda+2\mu)\int_{0}^{t}a(\tau)(a^\prime(\tau))^2 \d\tau \leq E_0,
\end{equation}
and
\begin{equation*}
a(t)\in H^1([0,T]),\ c_0\leq a(t) \leq C_T ,\ t\in(0,T),
\end{equation*}
where $E_0:=\int_{\varepsilon}^{a_0}(\frac{1}{2}\rho_0 u_0^2+\frac{1}{\gamma-1}\rho_0^\gamma) r^2\mathrm{d}r= \int_{0}^{\f{M}{4\pi}}(\frac{1}{2}u_0^2+\frac{1}{\gamma-1}\rho_0^{\gamma-1})\mathrm{d}x,$ $c_0$ and $C_T$ are positive constants.
 \end{lemma}
Then, we establish the uniform estimates of solutions away from symmetric center in Lagrangian coordinates.
\begin{lemma}\label{lemb1}
 Let $T>0$, $\f{6}{5}<\gamma\leq\f{4}{3}$ and $(\rho,u,a)$  be any strong solution to the FBVP \eqref{1.2a}-\eqref{3.2} for $\tau\in[0,T]$ under the assumption $M<M_c$.  Assume further that for $0<r_0<a_0,$
\begin{equation}\label{4.11}
    \rho_0(r)\in L^{\infty}([r_0,a_0]),\  u_0(r)\in H^1([r_0,a_0]),
\end{equation}
namely
\begin{equation*}
\rho_0(x)\in L^{\infty}([x_0,\f{M}{4\pi}]),\  (\f{u_0}{\sqrt{\rho_0}r(0)},\sqrt{\rho_0}r^2(0) u_{0,x})\in L^2([x_0,\f{M}{4\pi}]),
\end{equation*}
where $r^3(0):=r^3(x,0)=\varepsilon^3+3\int_{0}^{x}\frac{1}{\rho_0}(y)\d y.$
Then
  \begin{equation}\label{b1}
    \int_{0}^{T}\int_{x_1}^{\f{M}{4\pi}} (u^2_{\tau}+F^2_x)\d x\d\tau+\int_{x_1}^{\f{M}{4\pi}}\f{F^2}{\rho}\d x+\int_{x_1}^{\f{M}{4\pi}} (\rho r^4 |u_x|^2+ \f{2u^2}{\rho r^2}+|u_x|)\d x\leq C_{x_0}\bar{\delta}_{0},
  \end{equation}
where $ C_{x_0}$ is a positive constant just depending on T, $x_1$, $x_0 ~(0<x_0<x_1<\f{M}{4\pi})$, $E_0$ and $\bar{\delta}_{0}=\|\rho_0\|_{L^{\infty}[r_0,a_0]}+\|u_0\|_{H^1[r_0,a_0]}=\|\rho_0\|_{L^{\infty}[x_0,\f{M}{4\pi}]}
+\|(\f{u_0}{\sqrt{\rho_0}r(0)},\sqrt{\rho_0}r^2(0) u_{0,x})\|_{L^{2}[x_0,\f{M}{4\pi}]}$.
\end{lemma}
\proof Multiplying $(\ref{3.10})_2$ with $u_{\tau}\phi$ and integrating the resulted equation over $[0,\f{M}{4\pi}]$, then we have
\begin{equation}\label{b2}
   \int_{0}^{\f{M}{4\pi}} u^2_{\tau}\phi\d x+\int_{0}^{\f{M}{4\pi}} F_x (u r^2)_{\tau} \phi\d x-\int_{0}^{\f{M}{4\pi}} F_x 2ru^2 \phi\d x+\int_{0}^{\f{M}{4\pi}} \f{4\pi x}{r^2}u_\tau \phi\d x=0,
\end{equation}
where $\phi=\chi^2(x)$ and $\chi\in C^\infty([0,\f{M}{4\pi}])$ satisfies $0\leq\chi(x)\leq 1,\ \chi(x)=1$ for $x\in[x_1,\f{M}{4\pi}](0<x_0<x_1<\f{M}{4\pi}),\ \chi(x)=0$ for $x\in[0,x_0]$ and $|\chi^\prime|\leq \frac{2}{x_1-x_0}.$ After integrating by part, it holds that
\begin{eqnarray}\label{b3}
 \nonumber &&\int_{0}^{\f{M}{4\pi}} u^2_{\tau}\phi\d x+\int_{0}^{\f{M}{4\pi}} F (\f{F-\rho^{\gamma}}{(\lambda+2\mu)\rho})_{\tau} \phi\d x-\int_{0}^{\f{M}{4\pi}} F (u r^2)_{\tau} \phi^{\prime}\d x\\
&&+\int_{0}^{\f{M}{4\pi}} (r^{-2}u_{\tau} +4\pi xr^{-4})2ru^2 \phi\d x+\int_{0}^{\f{M}{4\pi}} \f{4\pi x}{r^2}u_\tau \phi\d x=0,
\end{eqnarray}
which implies
\begin{eqnarray}\label{b4}
  \nonumber&&\int_{0}^{\f{M}{4\pi}} u^2_{\tau}\phi\d x+\f{1}{2(\lambda+2\mu)}\f{\d}{\d \tau}\int_{0}^{\f{M}{4\pi}} \f{F^2}{\rho} \phi\d x \\
 \nonumber&&=\f{1}{2(\lambda+2\mu)}\int_{0}^{\f{M}{4\pi}}F \f{F}{\sqrt{\rho}}\f{\rho_{\tau}}{\rho^{\f{3}{2}}} \phi\d x+\f{\gamma-1}{\lambda+2\mu}\int_{0}^{\f{M}{4\pi}} F \rho^{\gamma-2}\rho_{\tau} \phi\d x-2\int_{0}^{\f{M}{4\pi}} r^{-1}u_{\tau} u^2 \phi\d x\\
&&+\int_{0}^{\f{M}{4\pi}} F u_{\tau}r^2 \phi^{\prime}\d x+2\int_{0}^{\f{M}{4\pi} }F u^2 r \phi^{\prime}\d x-8\pi\int_{0}^{\f{M}{4\pi}}u^2r^{-3}x\phi\d x-4\pi\int_{0}^{\f{M}{4\pi}}u_{\tau}r^{-2}x\phi\d x,
\end{eqnarray}
where we have used the fact
\begin{equation*}
    -(u r^2)_x=\f{F-\rho^{\gamma}}{(\lambda+2\mu)\rho},\ F_x=-r^{-2}u_{\tau}-4\pi xr^{-4}.
\end{equation*}
Using Lemma \ref{lem2.10}-\ref{lem2.3}, H\"{o}lder's  inequality, we obtain
\begin{eqnarray}\label{b5}
 \nonumber  &&\int_{0}^{\f{M}{4\pi}} u^2_{\tau}\phi\d x+\f{1}{2(\lambda+2\mu)}\f{\d}{\d \tau}\int_{0}^{\f{M}{4\pi}} \f{F^2}{\rho} \phi\d x\\
  \nonumber&&\leq\delta\int_{0}^{\f{M}{4\pi}} u^2_{\tau}\phi\d x +\delta\|\sqrt{\phi}F\|^2_{L^\infty[0,\f{M}{4\pi}]}+\delta^{-1}C(\int_{x_0}^{\f{M}{4\pi}} \rho^{-3}|\rho_{\tau}|^2 \d x)\int_{0}^{\f{M}{4\pi}} \f{F^2}{\rho} \phi\d x \\
    \nonumber&&+C\int_{x_0}^{\f{M}{4\pi}} F^2\d x+C_{x_0}\int_{x_0}^{\f{M}{4\pi}}\rho^{-3}|\rho_{\tau}|^2 +\delta^{-1}C_{x_0}\|u\sqrt{\phi}\|^2_{L^\infty[0,\f{M}{4\pi}]}\int_{0}^{\f{M}{4\pi}} u^2\d x,\\
 &&+\delta^{-1}C_{x_0}\int_{x_0}^{\f{M}{4\pi}} F^2\d x+C_{x_0}\|u\sqrt{\phi}\|^2_{L^\infty[0,\f{M}{4\pi}]}\int_{0}^{\f{M}{4\pi}} u^2\d x+ C_{x_0}\int_{0}^{\f{M}{4\pi}} u^2\d x+\delta^{-1}C_{x_0},
\end{eqnarray}
where $\delta\in(0,1)$ is a small positive constant and $C_{x_0}$ is a constant depending on $x_0,~ E_0,~T$ and $\|\rho_0\|_{L^{\infty}[x_0,\f{M}{4\pi}]}$.
Since
\begin{equation}\label{b6}
   \sqrt{\phi}F=-\int_{x}^{\f{M}{4\pi}}(\sqrt{\phi}F)_y \d y=\int_{x}^{\f{M}{4\pi}}u_{\tau}r^{-2}\sqrt{\phi}\d y
+4\pi\int_{x}^{\f{M}{4\pi}}y r^{-4}\sqrt{\phi}\d y-\int_{x}^{\f{M}{4\pi}}(\sqrt{\phi})^\prime F \d y,
\end{equation}
then we obtain
\begin{equation}\label{b7}
\|\sqrt{\phi}F\|^2_{L^\infty[0,\f{M}{4\pi}]}\leq C_{x_0}\int_{0}^{\f{M}{4\pi}} u^2_{\tau}\phi\d x+C_{x_0}\int_{x_0}^{\f{M}{4\pi}} F^2\d x+ C_{x_0}.
\end{equation}
By \eqref{3.19a}, $(\ref{b5})$ and $(\ref{b7})$, we obtain
\begin{eqnarray}\label{b8}
 \nonumber\int_{0}^{\f{M}{4\pi}} u^2_{\tau}\phi\d x+\f{\d}{\d \tau}\int_{0}^{\f{M}{4\pi}} \f{F^2}{\rho} \phi\d x &\leq& C_{x_0}(\int_{x_0}^{\f{M}{4\pi}} \rho^{-3}|\rho_{\tau}|^2 \d x)\int_{0}^{\f{M}{4\pi}} \f{F^2}{\rho} \phi\d x+C_{x_0}\|u\phi\|^2_{L^\infty[0,\f{M}{4\pi}]} \\
  &&+C_{x_0}\int_{x_0}^{\f{M}{4\pi}} (F^2+\rho^{-3}|\rho_{\tau}|^2)\d x+\delta^{-1}C_{x_0},
\end{eqnarray}
By $(\ref{3.10})_1$, Lemma \ref{lem2.10}-\ref{lem2.3}, we get
\begin{equation}\label{b9}
    \int_{0}^{T} \int_{x_0}^{\f{M}{4\pi}} F^2\d x\d\tau+ \int_{0}^{T}\int_{0}^{\f{M}{4\pi}} \rho^{-3}|\rho_{\tau}|^2\d x\d\tau \leq C_{x_0} .
\end{equation}
\begin{eqnarray}\label{jb9}
 \nonumber  &&\int_{0}^{T}\|u\sqrt{\phi}\|^2_{L^\infty[0,\f{M}{4\pi}]}\d\tau\leq C \int_{0}^{T}(\int_{0}^{\f{M}{4\pi}}|u|\d x+\int_{0}^{\f{M}{4\pi}}|u_x|\sqrt{\phi}\d x)^2 \d\tau \\
\nonumber&&\leq C \int_{0}^{T}\int_{0}^{\f{M}{4\pi}}u^2\d x\d\tau+C \int_{0}^{T}\int_{0}^{\f{M}{4\pi}}\rho r^4|u_x|^2\d x\d\tau+C_{x_0} \int_{0}^{T}\int_{0}^{\f{M}{4\pi}}\f{1}{\rho r^2}\phi\d x\d\tau\\
&&\leq C \int_{0}^{T}\int_{0}^{\f{M}{4\pi}}\rho r^4|u_x|^2\d x\d\tau+C_{x_0}\leq C_{x_0}.
\end{eqnarray}
It follows $(\ref{b8})$-$(\ref{jb9})$ and Gronwall's inequality that
\begin{equation}\label{b10}
   \int_{0}^{T}\int_{0}^{\f{M}{4\pi}} u^2_{\tau}\phi\d x\d\tau+\int_{0}^{\f{M}{4\pi}} \f{F^2}{\rho} \phi\d x\leq C_{x_0}\bar{\delta}_{0}.
\end{equation}
The fact that
\begin{eqnarray}\label{b11}
 \nonumber\int_{0}^{\f{M}{4\pi}} \f{|\rho(r^2u)_x|^2}{\rho} \phi\d x &=& \int_{0}^{\f{M}{4\pi}} (\rho r^4 |u_x|^2+ \f{2u^2}{\rho r^2})\phi\d x+2u^2(\f{M}{4\pi},\tau)a(\tau) \\
&&-2\int_{0}^{\f{M}{4\pi}}r u^2\phi^{\prime}\d x,
\end{eqnarray}
 together with $(\ref{b10})$ gives rise to
\begin{eqnarray}\label{b12}
 \nonumber &&\int_{0}^{\f{M}{4\pi}} (\rho r^4 |u_x|^2+ \f{2u^2}{\rho r^2})\phi\d x+u^2(\f{M}{4\pi},\tau)a(\tau)\\
 \nonumber&&\leq C\int_{0}^{\f{M}{4\pi}} \f{((\lambda+2\mu)\rho(r^2u)_x)^2}{\rho} \phi\d x+C\int_{0}^{\f{M}{4\pi}}u^2\d x\\
 \nonumber&&\leq C\int_{0}^{\f{M}{4\pi}} \f{F^2}{\rho} \phi\d x+C\int_{0}^{\f{M}{4\pi}} \rho^{2\gamma-1} \phi\d x+C\int_{0}^{\f{M}{4\pi}}u^2\d x\\
&&\leq C\int_{0}^{\f{M}{4\pi}} \f{F^2}{\rho} \phi\d x+ C_{x_0}\leq  C_{x_0}\bar{\delta}_{0}.
\end{eqnarray}
Therefore, it holds that
\begin{equation}\label{jb10}
  \int_{0}^{\f{M}{4\pi}}|u_x|\phi\d x\leq C\int_{0}^{\f{M}{4\pi}}\rho r^4 |u_x|^2\phi\d x+ C_{x_0}\int_{0}^{\f{M}{4\pi}}\f{1}{\rho}\phi\d x\leq  C_{x_0}\bar{\delta}_{0}.
\end{equation}
\begin{lemma}\label{lemb2}
Let $T>0$, $\f{6}{5}<\gamma\leq\f{4}{3}$, $(\rho,u,a)$ be the solution to FBVP \eqref{1.2a}-\eqref{3.2} for $(r,t)\in[\varepsilon,a(t)]\times[0,T]$ under the assumption $M<M_c$.
If the initial data $\rho_0$ also satisfies $(\rho_0^{q})_x\in L^2[x_1,\f{M}{4\pi}],$ then
 \begin{equation}\label{jb22}
  \int_{x_1}^{\f{M}{4\pi}}|(\rho^{q})_x|^2\d x\leq C_{x_1}\delta_0,
\end{equation}
where $\f{1}{2}<q=k+\f{1}{2}\leq\gamma$ and the constant $ C_{x_1}$ depends on T, $x_1$, $x_0 (0<x_0<x_1<\f{M}{4\pi})$, $E_0$ and  $\delta_0=\bar{\delta}_{0}+\|(\rho_0^{q})_x\|_{L^2[x_1,\f{M}{4\pi}]}$.
\end{lemma}
\proof Multiplying  $(\ref{3.10})_1$ by $q\rho^{q-1}$ and differentiating the resulted equation with respect to $x$, we obtain
\begin{equation}\label{b19}
    \rho^{q}_{x\tau}+q(\rho^{q+1}(ur^2)_x)_x=0.
\end{equation}
Multiplying $(\ref{b19})$ by $(\rho^{q})_x$ and integrating over $[x_1,\f{M}{4\pi}]$, we have
\begin{equation}\label{b20}
     \f{1}{2}\f{\d}{\d\tau}\int_{x_1}^{\f{M}{4\pi}}|(\rho^{q})_x|^2\d x=-q\int_{x_1}^{\f{M}{4\pi}}\rho(ur^2)_x|(\rho^q)_x|^2\d x-\f{q}{\lambda+2\mu}\int_{x_1}^{\f{M}{4\pi}}\rho^q(\rho^q)_x(r^{-2}u_{\tau}+(\rho^{\gamma})_x+4\pi xr^{-4})\d x,
\end{equation}
from which, together with Lemma \ref{lem2.2}-\ref{lem2.3}, Lemma \ref{lem2.10}-\ref{lemb1}, we get
\begin{eqnarray}\label{b21}
 \nonumber&&\f{1}{2}\f{\d}{\d\tau}\int_{x_1}^{\f{M}{4\pi}}|(\rho^{q})_x|^2\d x
\leq C\|\rho(ur^2)_x\|_{L^{\infty}[x_1,\f{M}{4\pi}]}\int_{x_1}^{\f{M}{4\pi}}|(\rho^{q})_x|^2\d x
+C\int_{x_1}^{\f{M}{4\pi}}\rho^{\gamma}|(\rho^{q})_x|^2\d x\\
 \nonumber&&+C\|r^{-2}\|_{L^{\infty}[x_1,\f{M}{4\pi}]}\int_{x_1}^{\f{M}{4\pi}}(|(\rho^{q})_x|^2+\rho^{2q}u_{\tau}^2)\d x
+C\|r^{-4}\|_{L^{\infty}[x_1,\f{M}{4\pi}]}\int_{x_1}^{\f{M}{4\pi}}(|(\rho^{q})_x|^2+\rho^{2q})\d x\\
\nonumber&&\leq C_{x_1}(\|F\|_{L^{\infty}[x_1,\f{M}{4\pi}]}+\|\rho\|^{\gamma}_{L^{\infty}[x_1,\f{M}{4\pi}]})\int_{x_1}^{\f{M}{4\pi}}|(\rho^{q})_x|^2\d x
+C\int_{x_1}^{\f{M}{4\pi}}\rho^{\gamma}|(\rho^{q})_x|^2\d x+C_{x_1} \int_{x_1}^{\f{M}{4\pi}}\rho^{2q}(u_{\tau}^2+1)\d x\\
  &&\leq C_{x_1}(\int_{x_1}^{\f{M}{4\pi}}|u_{\tau}|^2\d x+1)\int_{x_1}^{\f{M}{4\pi}}|(\rho^{q})_x|^2\d x
+C_{x_1}\int_{x_1}^{\f{M}{4\pi}}|u_{\tau}|^2\d x+C_{x_1},
\end{eqnarray}
where $C_{x_1}$ is a constant depending on $x_1,~ E_0,~T$ and $\|\rho_0\|_{L^{\infty}[x_0,\f{M}{4\pi}]}$.\\
By Gronwall's inequality and Lemma \ref{lemb1}, it holds that
\begin{equation}\label{b22}
  \int_{x_1}^{\f{M}{4\pi}}|(\rho^{q})_x|^2\d x\leq C_{x_1}\delta_0.
\end{equation}
\begin{lemma}\label{lemb3}
 Let $T>0$, $\f{6}{5}<\gamma\leq\f{4}{3}$ and $(\rho,u,a)$  be any strong solution to the FBVP \eqref{1.2a}-\eqref{3.2}for $(r,t)\in[\varepsilon,a(t)]\times[0,T]$  under the assumption $M<M_c$.
Assume further $\rho_0^{-\f{1}{2}}r^{-2}\partial_r^2u_{0}(r)\in L^2[r_1,a_0]$, namely $(\rho_0 r^2 u_{0x}(x))_x\in L^2[x_1,\f{M}{4\pi}]$ with
$x_1=\int_{\varepsilon}^{r_1}\rho_0 s^2\d s$,  then
\begin{equation}\label{b13}
\int_{x_2}^{\f{M}{4\pi}} u^2_{\tau}\d x+\int_{0}^{T}\int_{x_2}^{\f{M}{4\pi}}\f{F_{\tau}^2}{\rho}\d x\d\tau+\int_{0}^{T}\int_{x_2}^{\f{M}{4\pi}}(\rho r^4 u_{x\tau}^2+\f{2 u^2_{\tau}}{\rho r^2})\d x\d\tau+2\int_{0}^{T}|a^{\prime\prime}(\tau)|^2\d\tau\leq C_{x_2}\delta_1 ,
\end{equation}
and
\begin{equation}\label{jb13}
   \int_{x_2 }^{\f{M}{4\pi}}|(\rho (u r^2)_{x})_{x}|^2\d x\leq C_{x_2}\delta_1,
\end{equation}
where $C_{x_2}$  depends on T, $x_0$, $x_1$, $x_2$ $(0<x_0<x_1<x_2<\f{M}{4\pi})$, $E_0$ and $\delta_1=\delta_0
+\|\rho_0 r^2 u_{0x}\|_{L^2[x_1,\f{M}{4\pi}]}$.
\end{lemma}
\proof Differentiating $(\ref{3.10})_2$ with respect to $\tau$, we obtain
\begin{equation}\label{b14}
  u_{\tau\tau}+r^2F_{x\tau}+2ruF_x-8\pi xr^{-3}u=0.
\end{equation}
Choose a smooth function  $\psi=\zeta^2(x)$ where $\zeta\in C^\infty([0,\f{M}{4\pi}])$ satisfies $0\leq\zeta(x)\leq 1,\ \zeta(x)=1$ for $x\in[x_2,\f{M}{4\pi}](0<x_0<x_1<x_2<\f{M}{4\pi}),\ \zeta(x)=0$ for $x\in[0,x_1]$ and $|\zeta^\prime|\leq \frac{2}{x_2-x_1}.$\\
Taking inner product of $(\ref{b14})$ with $u_{\tau}\psi$ and integrating by part, it holds that
\begin{eqnarray}\label{b15}
\nonumber&& \f{1}{2}\f{\d}{\d\tau}\int_{0}^{\f{M}{4\pi}} u^2_{\tau}\psi\d x+\f{1}{\lambda+2\mu}\int_{0}^{\f{M}{4\pi}}\f{F_{\tau}^2}{\rho}\psi\d x\\
\nonumber& =&\f{1}{\lambda+2\mu}\int_{0}^{\f{M}{4\pi}}\f{FF_{\tau}}{\rho^2}\rho_{\tau}\psi\d x+\int_{0}^{\f{M}{4\pi}}F_{\tau}u_{\tau}r^2\psi^{\prime}\d x\d\tau-\f{4}{\lambda+2\mu}\int_{0}^{\f{M}{4\pi}}F_{\tau}\rho^{\gamma}\f{u}{\rho r}\psi\d x\\
\nonumber&&+\f{\gamma-1}{(\lambda+2\mu)^2}\int_{0}^{\f{M}{4\pi}}F_{\tau}F\rho^{\gamma-1}\psi\d x-\f{\gamma-1}{(\lambda+2\mu)^2}\int_{0}^{\f{M}{4\pi}}F_{\tau}\rho^{2\gamma-1}\psi\d x+\f{4}{\lambda+2\mu}\int_{0}^{\f{M}{4\pi}}F_{\tau}F\f{u}{\rho r}\psi\d x\\
\nonumber&&+6\int_{0}^{\f{M}{4\pi}}F_{\tau}\f{u^2}{\rho r^2}\psi\d x+2\int_{0}^{\f{M}{4\pi}} u^2_{\tau}ur^{-1}\psi\d x+16\pi\int_{0}^{\f{M}{4\pi}}x u u_{\tau}r^{-3}\psi\d x\\
\nonumber&\leq&\delta\int_{0}^{\f{M}{4\pi}}\f{F_{\tau}^2}{\rho}\psi\d x+C(\|u\|_{L^\infty([x_1,\f{M}{4\pi}]}+1)\int_{0}^{\f{M}{4\pi}} u^2_{\tau}\psi\d x+\delta^{-1}C\|F\|^2_{L^\infty([x_1,\f{M}{4\pi}]}\int_{x_1}^{\f{M}{4\pi}}(\rho^{-3}|\rho_{\tau}|^2+\f{u^2}{\rho r^2})\d x\\
\nonumber&&+\delta^{-1}C(\|\rho\|^{2\gamma}_{L^\infty([x_1,\f{M}{4\pi}]}+\|u\|^2_{L^\infty([x_1,\f{M}{4\pi}]})\int_{x_1}^{\f{M}{4\pi}}\f{u^2}{\rho r^2}\d x+\delta^{-1}C\|\rho\|_{L^\infty([x_1,\f{M}{4\pi}]}\int_{x_1}^{\f{M}{4\pi}}u^2_{\tau}\d x\\
\nonumber&&+\delta^{-1}C\|\rho\|^{3\gamma}_{L^\infty([x_1,\f{M}{4\pi}]}\int_{0}^{\f{M}{4\pi}}\rho^{\gamma-1}\psi\d x+\delta^{-1}C\|\rho\|^{2\gamma}_{L^\infty([x_1,\f{M}{4\pi}]}\int_{0}^{\f{M}{4\pi}}F^2\psi\d x+C\int_{x_1}^{\f{M}{4\pi}}u^2\d x\\
&\leq& \delta\int_{0}^{\f{M}{4\pi}}\f{F_{\tau}^2}{\rho}\psi\d x+C(\|u\|_{L^\infty([x_1,\f{M}{4\pi}]}+1)\int_{0}^{\f{M}{4\pi}} u^2_{\tau}\psi\d x+
\delta ^{-1}C_{x_2}(\|u\|^2_{L^\infty([x_1,\f{M}{4\pi}]}+\int_{x_1}^{\f{M}{4\pi}}u^2_{\tau}\d x+1),
\end{eqnarray}
where we have used the following facts:
\begin{equation*}
(ur^2)_x=\f{1}{\lambda+2\mu}\cdot\f{\rho^\gamma-F}{\rho},~~~~(ur^2)_{x\tau}=\f{1}{\lambda+2\mu}\cdot(\f{\rho^\gamma-F}{\rho})_{\tau},
\end{equation*}
\begin{equation}\label{b16}
    \int_{x_0}^{\f{M}{4\pi}}\rho^{-3}|\rho_{\tau}|^2\d x\leq\int_{0}^{\f{M}{4\pi}} (\rho r^4 |u_x|^2+ \f{2u^2}{\rho r^2})\phi\d x\leq C,
\end{equation}
\begin{eqnarray}\label{b17}
  \nonumber|F(y,\tau)|&=&|-\int_{y}^{\f{M}{4\pi}}F_x\d x|=|\int_{y}^{\f{M}{4\pi}}r^{-2}u_{\tau}\d x+4\pi\int_{y}^{\f{M}{4\pi}}r^{-4}x\d x|\\
  &\leq &\f{C}{r^2(x_1,\tau)}(\int_{x_1}^{\f{M}{4\pi}}u^2_{\tau}\d x)^{\f{1}{2}}+\f{C}{r^4(x_1,\tau)},\ y\in[x_1,\f{M}{4\pi}].
\end{eqnarray}
Choosing a small $\delta\in(0,1)$ and using Gronwall's inequality, we get
\begin{equation}\label{b18}
  \int_{0}^{\f{M}{4\pi}} u^2_{\tau}\psi\d x+\int_{0}^{T}\int_{0}^{\f{M}{4\pi}}\f{F_{\tau}^2}{\rho}\psi\d x\d\tau\leq C_{x_2}\delta_1.
\end{equation}
Furthermore, we have
\begin{eqnarray}\label{b19a}
\nonumber&&\int_{0}^{T}\int_{0}^{\f{M}{4\pi}}(\rho r^4 u_{x\tau}^2+\f{2 u^2_{\tau}}{\rho r^2})\psi\d x\d\tau+2\int_{0}^{T}|a^{\prime\prime}(\tau)|^2a(\tau)\d\tau\\
 \nonumber &=&\int_{0}^{T}\int_{0}^{\f{M}{4\pi}}\f{(\rho r^2 u_{x\tau}+\f{2 u_{\tau}}{r})^2}{\rho}\psi\d x\d\tau+2\int_{0}^{T}\int_{0}^{\f{M}{4\pi}}u_{\tau}^{2}r\psi^{\prime}\d x\d\tau\\
 \nonumber &\leq&\int_{0}^{T}\int_{0}^{\f{M}{4\pi}}\f{|(\rho r^2 u_{x}+\f{2 u}{r})_{\tau}|^2}{\rho}\psi\d x\d\tau+\int_{0}^{T}\int_{0}^{\f{M}{4\pi}}\f{|\rho_{\tau} r^2 u_{x}|^2}{\rho}\psi\d x\d\tau+4\int_{0}^{T}\int_{0}^{\f{M}{4\pi}}\f{|\rho r u u_{x}|^2}{\rho}\psi\d x\d\tau\\
 \nonumber &&+C\int_{0}^{T}\int_{0}^{\f{M}{4\pi}}\f{u^{4}}{\rho r^2}\psi\d x\d\tau+2\int_{0}^{T}\int_{0}^{\f{M}{4\pi}}u_{\tau}^{2}r|\psi^{\prime}|\d x\d\tau\\
 \nonumber &\leq&\int_{0}^{T}\int_{0}^{\f{M}{4\pi}}\f{|F_{\tau}|^2}{\rho}\psi\d x\d\tau+\int_{0}^{T}\int_{0}^{\f{M}{4\pi}}\rho^{2\gamma-3}|\rho_{\tau}|^2\psi\d x\d\tau+\int_{0}^{T}\|\rho r^2 u_x\|^2_{L^\infty([x_1,\f{M}{4\pi}]}\int_{x_0}^{\f{M}{4\pi}}\rho^{-3}|\rho_{\tau}|^2\d x\d\tau\\
&&+\int_{0}^{T}(\|\rho r^2 u_x\|^2_{L^\infty([x_1,\f{M}{4\pi}]}+\|u\|^2_{L^\infty([x_1,\f{M}{4\pi}]})\int_{x_1}^{\f{M}{4\pi}}\f{u^2}{\rho r^2}\d x\d\tau+2\int_{0}^{T}\int_{x_1}^{\f{M}{4\pi}}u_{\tau}^{2}\d x\d\tau\leq C_{x_2}\delta_1.
\end{eqnarray}
By $(\ref{b18})$ and Lemma $\ref{lemb2},$ we obtain
\begin{eqnarray}\label{B20}
  \nonumber&&\int_{0}^{\f{M}{4\pi}}|(\rho (u r^2)_{x})_{x}|^2\psi\d x\leq C\int_{0}^{\f{M}{4\pi}}|F_{x}|^2\psi\d x+C\int_{0}^{\f{M}{4\pi}}|(\rho^{\gamma})_{x}|^2\psi\d x\\
 &&\leq C_{x_2}\int_{0}^{\f{M}{4\pi}}|u_{\tau}|^2\psi\d x+C_{x_2}\int_{0}^{\f{M}{4\pi}}|\f{4\pi x}{r^4}|^2\psi\d x+C_{x_2} \int_{x_1}^{\f{M}{4\pi}}|(\rho^{\gamma})_{x}|^2\d x\leq C_{x_2}\delta_1.
\end{eqnarray}
At last, we obtain the interior estimates below.
\begin{lemma}\label{lemb4}
 Let $T>0$, $\f{6}{5}<\gamma\leq\f{4}{3}$ and $(\rho,u,a)$  be any strong solution to the FBVP \eqref{1.2a}-\eqref{3.2}for $(r,t)\in[\varepsilon,a(t)]\times[0,T]$ under the assumption $M<M_c$. Assume further there exists $0<r_0^-<r_0<r_b<r_b^+\leq a_0$  and a positive constant $\rho_{*}$ such that
\begin{equation}\label{in}
    \inf_{r\in[r_0^-,r_b^+]}\rho_0(r)\geq\rho_{*}>0,\ u_0(r)\in H^2([r_0^-,r_b^+]),
\end{equation}
that is
\begin{equation}\label{in0}
   \inf_{x\in[x_0^-,x_b^+]}\rho_0(x)\geq\rho_{*}>0,\ u_0(x)\in H^2([x_0^-,x_b^+]),
\end{equation}
then it holds that
\begin{equation}\label{in1}
  0< c_{x_0^-,T}\leq\rho(r,t)\leq C_{x_0^-,T},\ \forall t\in [0,T],\ \forall r\in [r_{x_0^-}(t),r_{x_b^+}(t)],
\end{equation}
  \begin{eqnarray}\label{in2}
  \nonumber&&\sup_{\tau\in[0,T]}(\|u_x\|_{L^2[x_0,x_b]}^2+\|\rho_x\|_{L^2[x_0,x_b]}^2+\|\rho_\tau\|_{L^2[x_0,x_b]}^2+\|F\|_{L^2[x_0,x_b]}^2)  \\
  && +\int_{0}^{T}\|(u_{\tau},F_x)(\tau)\|_{L^2[x_0,x_b]}^{2}\mathrm{d}\tau+\int_{0}^{T}\|(u_{xx},\rho_{x\tau}(\tau)\|_{L^2[x_0,x_b]}^{2}\mathrm{d}\tau\leq C_1,
  \end{eqnarray}
 \begin{eqnarray}\label{in3}
  \nonumber && \sup_{\tau\in[0,T]}(\|u_\tau\|_{L^2[x_0,x_b]}^2+\|F_x\|_{L^2[x_0,x_b]}^2+\|u_{xx}\|_{L^2[x_0,x_b]}^2)+\int_{0}^{T}\|F_{\tau}(\tau)\|_{L^2[x_0,x_b]}^{2}\mathrm{d}\tau\\
&&+\int_{0}^{T}\|(u_{\tau x},F_{xx})(\tau)\|_{L^2[x_0,x_b]}^{2}\mathrm{d}\tau\leq C_2,
 \end{eqnarray}
where $c_{x_0^-,T},\ C_{x_0^-,T}$ are positive constants depending on $x_0^-,\ T$ and the initial data, and $r_{x_0^-}(t),\ r_{x_b^+}(t)$ are particle paths with $x_0^-=\f{M}{4\pi}-\int_{r_0^-}^{a_0}\rho_0 r^2\d r,\ x_b^+=\f{M}{4\pi}-\int_{r_b^+}^{a_0}\rho_0 r^2\d r$ and $x_i=\f{M}{4\pi}-\int_{r_i}^{a_0}\rho_0 r^2\d r\ (i=0,b).$
The constant $C_1$ depends on $E_0,\ c_{x_0^-,T},\ \ C_{x_0^-,T},\ x_b^+,\ x_i(i=0,b),\ \|\rho_0\|_{H^1[x_0^-,x_b^+]}$ and $\|u_0\|_{H^1[x_0^-,x_b^+]}$ and $C_2$ depends on $C_1$ and $\|u_0\|_{H^2[x_0^-,x_b^+]}$.
\end{lemma}
\proof \eqref{in1} can be deduced from \eqref{3.25} and \eqref{in}.\\
Multiplying $(\ref{3.10})_2$ with $u_{\tau}\widetilde{\phi}$ and integrating the resulted equation over $[0,\f{M}{4\pi}]$, then we have
\begin{equation}\label{in4}
  \int_{0}^{\f{M}{4\pi}} u^2_{\tau}\widetilde{\phi}\d x+\int_{0}^{\f{M}{4\pi}} F_x (u r^2)_{\tau}\widetilde{\phi}\d x-\int_{0}^{\f{M}{4\pi}} F_x 2ru^2 \widetilde{\phi}\d x+\int_{0}^{\f{M}{4\pi}} \f{4\pi x}{r^2}u_\tau \widetilde{\phi}\d x=0,
\end{equation}
where $\widetilde{\phi}=\widetilde{\chi}^2(x)$ and $\widetilde{\chi}\in C^\infty([0,\f{M}{4\pi}])$ satisfies $0\leq\widetilde{\chi}(x)\leq 1,\ \widetilde{\chi}(x)=1$ for $x\in[\f{x_0+x_0^-}{2},\f{x_b+x_b^+}{2}],\ \widetilde{\chi}(x)=0$ for $x\in[0,x_0^-]\cup[x_b^+,\f{M}{4\pi}]$ and $|\widetilde{\chi}^\prime|\leq \frac{4}{x_0-x_0^-}+\frac{4}{x_b^+-x_b}.$
\ Similarly as the proof of Lemma \ref{lemb1}, we obtain
\begin{equation}\label{in5}
   \int_{0}^{T}\int_{0}^{\f{M}{4\pi}} u^2_{\tau}\widetilde{\phi}\d x\d\tau+\int_{0}^{\f{M}{4\pi}} \f{F^2}{\rho}\widetilde{\phi}\d x\leq C_0,
\end{equation}
and
\begin{equation}\label{in6}
   (\lambda+2\mu)^2\int_{0}^{\f{M}{4\pi}} (\rho r^4 |u_x|^2+ \f{2u^2}{\rho r^2})\widetilde{\phi}\d x=\int_{0}^{\f{M}{4\pi}} \f{((\lambda+2\mu)\rho(r^2u)_x)^2}{\rho} \widetilde{\phi}\d x+2(\lambda+2\mu)^2\int_{0}^{\f{M}{4\pi}}u^2\widetilde{\phi}^{\prime}\d x\leq C_0,
\end{equation}
where $C_0$ depends on $E_0,\ C_{x_0^-,T}$ and $\|u_0\|_{H^1[x_0^-,x_b^+]}$.\\
From \eqref{in5}-\eqref{in6}, $(\ref{3.10})_2$ and \eqref{in1}, we get
\begin{equation}\label{in7}
    \sup_{\tau\in[0,T]}(\|u_x\|_{L^2[\f{x_0+x_0^-}{2},\f{x_b+x_b^+}{2}]}^2+\|F\|_{L^2[\f{x_0+x_0^-}{2},\f{x_b+x_b^+}{2}]}^2)
    +\int_{0}^{T}\|(u_{\tau},F_x)(\tau)\|_{L^2[\f{x_0+x_0^-}{2},\f{x_b+x_b^+}{2}]}^{2}\mathrm{d}\tau\leq \widetilde{C}_0,
\end{equation}
 where $\widetilde{C}_0$ depends on $C_{0}$ and $c_{x_0^-,T}$.
Similarly as the Lemma \ref{lemb2}, we obtain
\begin{equation}\label{in8}
    \sup_{\tau\in[0,T]}\|\rho_x\|_{L^2[\f{x_0+x_0^-}{2},\f{x_b+x_b^+}{2}]}^2\leq C_1
\end{equation}
Thus, \eqref{in2} can be deduced from $(\ref{3.10})$ and \eqref{in7}-\eqref{in8}.
Similarly as the proof of the Lemma \ref{lemb3}, one can obtain \eqref{in3} with the help of a proper cut-off function and \eqref{in1}.\\
\subsection{ The proof of Theorem \ref{thm1}}
\underline{\it The proof of Theorem \ref{thm1}.}   For any fixed $\varepsilon$, the approximate FBVP problem (\ref{1.2a}) is essentially an one-dimensional problem.
Indeed, one can construct global solutions to the approximate FBVP (\ref{1.2a})-(\ref{3.2}). Using the uniform estimates established above, we can obtain the global solution to the original FBVP problem after compactness argument when $\varepsilon\rightarrow 0$ and justify the expected properties in Theorem \ref{thm1} for the limiting solution.\par
We can modify the initial data  $(\rho_0,u_0)$ in Theorem \ref{thm1} properly such that the modified initial data $(\rho^{\varepsilon}_{0},u^{\varepsilon}_{0})$  satisfies the following properties on $[\varepsilon,a_0]$:
\begin{equation}\label{5.1}
  \inf_{r\in[\varepsilon,a_0]}\rho^{\varepsilon}_{0}(r)>0
,\ \ u^{\varepsilon}_{0}(\varepsilon)=0,\ \ ((\rho^\varepsilon_{0})^\gamma-(\lambda+2\mu) u^\varepsilon_{0,r}-(\lambda+2\mu) \frac{2u^\varepsilon_{0}}{r})(a_{0})=0,
\end{equation}
 $((\rho^{\varepsilon}_{0})^{k},u^{\varepsilon}_{0})\rightarrow(\rho^k_0,u_0)$ strongly in $H^{1}([\varepsilon,a_0])$ as $\varepsilon\rightarrow 0_{+}$ and $\rho^{\varepsilon}_{0}(r)\to \rho_0(r)$ as $\varepsilon\to 0_+$, refer to \cite{SG} for construction of such function. One can apply the standard argument to obtain a unique strong solution local in time  and then by the a-priori estimates and a continuity argument, we can continue the local solution globally in time. Thus, the approximate FBVP (\ref{1.2a})-(\ref{3.2}) has a global strong solution $(\rho^\varepsilon,u^\varepsilon,a^\varepsilon)$ on the domain $[\varepsilon,a^\varepsilon(t)]\times[0,T]$ with the initial data $(\rho^{\varepsilon}_{0},u^{\varepsilon}_{0}).$  In addition, one also can use the similar space-discrete difference method \cite{SG} to prove the global existence.  Extend $(\rho^\varepsilon,u^\varepsilon)$ by setting $(\rho^\varepsilon,u^\varepsilon)(r,t)=(\rho^\varepsilon(\varepsilon,t),0)$ for $0\leq r\leq\varepsilon$ and denote the extension function still by $(\rho^\varepsilon,u^\varepsilon)$ for convenience.\par
First, we prove the strong convergence of $(\rho^\varepsilon,u^\varepsilon,a^\varepsilon)$ near the free boundary. It's enough to prove the strong convergence on the domain $[r_{x_b}^{\varepsilon},a^\varepsilon]\times[0,T]$, where $r=r_{x_b}^{\varepsilon}$ is a particle path with $r_{x_b}^{\varepsilon}(0)=r_b \in(b_0,a_0]$ and $x_b=\int_{r_b}^{a_0}\rho_0 r^2\mathrm{d}r$ and the initial data satisfies $(\rho^k_0,u_0)\in H^{1}[b_0,a_0].$
 It's convenient to show the strong convergence in Lagrangian coordinate on $[x_b,1]\times[0,T].$ Indeed, we can show $(\rho^\varepsilon,u^\varepsilon,a^\varepsilon)$ satisfies the uniform estimate established in Lemmas~\ref{lem2.10}--\ref{lemb3} on $[x_b,1]\times[0,T]$. Thus, by Lions-Aubin's lemma, there is a limiting function $(\rho_b,u_b,a)$ so that up to a subsequence $(\rho^{\varepsilon_j},u^{\varepsilon_j},a^{\varepsilon_j})$, it holds that
\begin{equation}\label{5.2}
   \left\{
     \begin{array}{ll}
        (\rho^{\varepsilon_j},u^{\varepsilon_j})\rightarrow(\rho_b,u_b)& \hbox{strongly in $C([0,T]\times[x_b,1])\times C([0,T];L^p[x_b,1]),$} \\
         F^{\varepsilon_j}\rightarrow F & \hbox{strongly in $ L^2([0,T]\times[x_b,1]),$} \\
       a^{\varepsilon_j}\rightarrow a & \hbox{strongly in $C^\alpha([0,1]),\alpha\in(0,\frac{1}{2}),$}
     \end{array}
   \right.
\end{equation}
where $r_\tau =u_b$ and $(r^3)_x=\frac{3}{\rho_b},\ F=\rho^{\gamma}_{b}-(\lambda+2\mu)\rho_b(r^2 u_b)_x=\rho^{\gamma}_{b}
-(\lambda+2\mu)\rho_br^2\partial_x u_b-(\lambda+2\mu)\frac{2u_b}{r}.$ In addition, by Lemma~\ref{lem2.3} and the construction that $\rho_0^\varepsilon(r)\to \rho_0(r)$ as $\varepsilon\to0_+$, we conclude that the boundary condition $\rho_b(a(t),t)=0$ holds.

Next, we show the convergence of$(\rho^{\varepsilon_j},u^{\varepsilon_j},a^{\varepsilon_j})$
on an interior domain $\Omega_{in}^{\varepsilon_j}$ defined by
$$\Omega_{in}^{\varepsilon_j}=:\{(r,t)|0\leq r< a^{\varepsilon_j}(t),\ 0\leq t\leq T\}\cap\{(r,t)|0\leq r\leq a(t),\ 0\leq t\leq T\}.$$
Due to the strong convergence (\ref{5.2}) of velocity and the particle path as $\varepsilon_j\rightarrow 0_+$, it holds that for $\varepsilon_j>0$ small enough
\begin{equation}\label{5.3}
    \Omega_{in}=:\{(r,t)|0\leq r\leq r_{x_{in}}(t),\ 0\leq t\leq T\}\subset\subset\Omega_{in}^{\varepsilon_j},
\end{equation}
where $r=r_{x_{in}}(t)$ is a particle path defined by
\begin{equation}\label{5.4}
    \frac{\mathrm{d}}{\mathrm{d}t}r_{x_{in}}(t)=u_b(r_{x_{in}}(t),t),\ r_{x_{in}}(0)=r_{in}\in (r_b,a_0),
\end{equation}
which satisfies that for $x_b<x_{in}=1-\int_{r_{in}}^{a_0}\rho_0 r^2\mathrm{d}r,$
\begin{equation}\label{5.5}
    0<c(x_{in}-x_{b})^{\frac{\gamma}{\gamma-1}}\leq r_{x_{in}}^3(t)- r_{x_b}^3(t),\ t\in[0,T].
\end{equation}
With help of Lemma \ref{lem4.1}, a proper cut-off function and a similar compactness argument as \cite{JZ2001}, we can show that there is a limiting function $(\rho_{in},u_{in})(r,t)\ ((r,t)\in \Omega_{in}),$ so that up to a sub-subsequence $(\rho^{\varepsilon_j},\rho^{\varepsilon_j}u^{\varepsilon_j})$ converge to
$(\rho_{in},\rho_{in}u_{in})$ in the sense that
\begin{equation}\label{5.6}
   \left\{
     \begin{array}{ll}
        \rho^{\varepsilon_j}\rightarrow \rho_{in}& \hbox{strongly in $L^p(0,T;\mathcal{L}^p(0,r_{in}(t))),\ \forall\ 1\leq p\leq 2\gamma$} \\
         \rho^{\varepsilon_j}u^{\varepsilon_j}\rightharpoonup \rho_{in}u_{in} & \hbox{weakly in $L^\infty(0,T;\mathcal{L}^{\frac{2\gamma}{\gamma+1}}(0,r_{in}(t))),$}
     \end{array}
   \right.
\end{equation}
and $(\rho_{in},u_{in})$ satisfies \eqref{1.2} on $\Omega_{in}$ in the sense of distribution. As \cite{JZ2001}, we define $\mathcal{L}^p(\Omega)):=\{f\in L^1_{loc}(\Omega)|\int_{\Omega}|f(r)|^p r^2\mathrm{d}r<\infty\}$ with norm $\|\cdot\|_{\mathcal{L}^p(\Omega))}:=(\int_{\Omega}|\cdot|^p r^2\mathrm{d}r)^{\frac{1}{p}}.$
Finally, define
\begin{equation}\label{5.7}
(\rho,\rho\mathbf{u})=
\left\{
  \begin{array}{ll}
    (\rho_b,\rho_b\mathbf{u}_b)(\mathbf{x},t), &\ r_{x_{b}}(t)\leq|\mathbf{x}|\leq a(t),\ t\in[0,T]\hbox{,} \\
    (\rho_{in},\rho_{in}\mathbf{u}_{in})(\mathbf{x},t), &\ 0\leq|\mathbf{x}|\leq r_{x_{in}}(t),\ t\in[0,T]\hbox{,}
  \end{array}
\right.
\end{equation}
where $\mathbf{u}=u\frac{\mathbf{x}}{r}$,\ $\mathbf{u_b}=u_b\frac{\mathbf{x}}{r}$,\ $\mathbf{u}_{in}=u_{in}\frac{\mathbf{x}}{r}$ and $r=|\mathbf{x}|.$  This is well defined and
\begin{equation}\label{5.8}
    (\rho_b,\rho_b\mathbf{u}_b)=(\rho_{in},\rho_{in}\mathbf{u}_{in}),\ a.e.\ (r,t)\in[r_{x_{b}(t)},r_{x_{in}}(t)]\times[0,T].
\end{equation}
We can easily deduce that $(\rho,\rho\mathbf{u},a(t))$ is a weak solution to FBVP \eqref{1.1}, \eqref{1.2}-\eqref{1.5} in the sense of Definition \ref{def},  and by similar argument to \cite{GLX2012} verify that $(\rho,\rho\mathbf{u},a(t))$ satisfies the properties (\ref{enegy1})-(\ref{2.15a}) and the free boundary condition with the help of Lemmas~\ref{lemb1}--\ref{lemb4}. The proof of Theorem \ref{thm1} is completed.

\section{Long Time Expanding Rate}
In this section, we investigate the long time behavior of global spherical symmetric solution to  the FBVP (\ref{1.2})-(\ref{1.5}). Indeed, we can obtain an expanding rate of the domain occupied by the fluid.\\

\underline{\it The proof of Theorem \ref{thm2}.}
Define an energy functional for a spherically symmetric solution $(\rho,u,a)$ as
\begin{eqnarray}\label{6.1}
  \nonumber H(t) &=& \int_{0}^{a(t)}(r-(1+t) u)^2\rho r^{2}\mathrm{d}r+\frac{2}{\gamma-1}(1+t) ^2\int_{0}^{a(t)}\rho^{\gamma} r^{2}\mathrm{d}r -(1+t) ^2\int_{0}^{a(t)}\f{4\pi}{r^{2}}(\int_{0}^{r}\rho s^{2}\d s)^2\mathrm{d}r\\
  \nonumber &=&\int_{0}^{a(t)}\rho r^{4}\mathrm{d}r-2(1+t) \int_{0}^{a(t)}\rho u r^{3}\mathrm{d}r+(1+t)^2\int_{0}^{a(t)} \rho u^2 r^{2}\mathrm{d}r\\
  &&+\frac{2}{\gamma-1}(1+t)^2\int_{0}^{a(t)}\rho^{\gamma} r^2\mathrm{d}r-(1+t) ^2\int_{0}^{a(t)}\f{4\pi}{r^{2}}(\int_{0}^{r}\rho s^{2}\d s)^2\mathrm{d}r.
\end{eqnarray}
A direct computation gives
\begin{eqnarray}\label{6.2}
 \nonumber &&H^\prime(t)= \int_{0}^{a(t)}(\rho_t r^{4}-2\rho u r^{3})\mathrm{d}r+ 2(1+t)\int_{0}^{a(t)} \{\rho u^2 r^2-(\rho u)_t r^3+\frac{2}{\gamma-1}\rho^{\gamma} r^2-\f{4\pi}{r^{2}}(\int_{0}^{r}\rho s^{2}\d s)^2\}\mathrm{d}r\\
 \nonumber &&+(1+t)^2\int_{0}^{a(t)} \{(\rho u^2)_t r^2+\frac{2}{\gamma-1}(\rho^{\gamma})_t r^2-\f{8\pi}{r^{2}}\int_{0}^{r}\rho s^{2}\d s\int_{0}^{r}\rho_t s^{2}\d s\}\mathrm{d}r\\
  \nonumber&&+\{\rho u r^{4}-2(1+t)\rho u^2 r^3+(1+t)^2\rho u^3 r^2+\frac{2}{\gamma-1}(1+t)^2\rho^{\gamma}u r^2-(1+t)^2\f{4\pi u}{r^{2}}(\int_{0}^{r}\rho s^{2}\d s)^2\}|_{r=a(t)}\\
  &&=:I_1+I_2+I_3+I_{B}.
\end{eqnarray}
By (\ref{1.2}) and (\ref{1.5}), one has
\begin{equation}\label{6.3}
   I_1 =-\int_{0}^{a(t)}((\rho u r^2)_r r^2+\rho u r^2 2r)\mathrm{d}r=-\int_{0}^{a(t)}(\rho u r^{4})_r\mathrm{d}r \\
  =-(\rho u r^{4})|_{r=a(t)},
\end{equation}

\begin{eqnarray}\label{6.5}
\nonumber  I_2 &=&2(1+t)\int_{0}^{a(t)} \{(\rho u^2 r^3)_r +(\rho^\gamma-(\lambda+2\mu) (u_r+\frac{2u}{r}))_r r^3+\frac{2}{\gamma-1}\rho^{\gamma} r^2\}\mathrm{d}r \\
\nonumber &&+2(1+t)\int_{0}^{a(t)} \{4\pi\rho r\int_{0}^{r}\rho s^{2}\d s-\f{4\pi}{r^{2}}(\int_{0}^{r}\rho s^{2}\d s)^2\}\mathrm{d}r  \\
  \nonumber &=& 2(1+t)(\rho u^2 r^3)|_{r=a(t)}-6(1+t)\int_{0}^{a(t)} (\rho^\gamma-(\lambda+2\mu) (u_r+\frac{2u}{r}))r^2\mathrm{d}r+\frac{4(1+t)}{\gamma-1}\int_{0}^{a(t)}\rho^{\gamma} r^2\mathrm{d}r\\
  \nonumber&&+2(1+t)\int_{0}^{a(t)} \{4\pi\rho r\int_{0}^{r}\rho s^{2}\d s-\f{4\pi}{r^{2}}(\int_{0}^{r}\rho s^{2}\d s)^2\}\mathrm{d}r\\
  \nonumber &=&6(\lambda+2\mu)(1+t)\int_{0}^{a(t)} (u r^2)_r\mathrm{d}r+\frac{2(2-3(\gamma-1))}{\gamma-1}(1+t)\int_{0}^{a(t)}\rho^{\gamma} r^2\mathrm{d}r\\
&&+2(1+t)\int_{0}^{a(t)} \{4\pi\rho r\int_{0}^{r}\rho s^{2}\d s-\f{4\pi}{r^{2}}(\int_{0}^{r}\rho s^{2}\d s)^2\}\mathrm{d}r+ 2(1+t)(\rho u^2 r^3)|_{r=a(t)},
\end{eqnarray}
\begin{eqnarray}\label{6.4}
  \nonumber I_3&=&(1+t)^2\int_{0}^{a(t)} 2u\rho u_t r^2 \mathrm{d}r+(1+t)^2\int_{0}^{a(t)} \rho_t u^2 r^2\mathrm{d}r\\
  \nonumber && +(1+t)^2\int_{0}^{a(t)} \frac{2\gamma}{\gamma-1}\rho^{\gamma-1}\rho_t r^2\mathrm{d}r-(1+t)^2\int_{0}^{a(t)} \f{8\pi}{r^{2}}\int_{0}^{r}\rho s^{2}\d s\int_{0}^{r}\rho_t s^{2}\d s\mathrm{d}r\\
 \nonumber&=&-(1+t)^2\int_{0}^{a(t)} 2u r^2\{\rho u u_r+\partial_r\rho^{\gamma}-(\lambda+2\mu)(u_r+\frac{2u}{r})_r+\f{4\pi\rho}{r^2}\int_{0}^{r}\rho s^2\d s\} \mathrm{d}r\\
 \nonumber&&-(1+t)^2\int_{0}^{a(t)} u^2(\rho u r^2)_r\mathrm{d}r-(1+t)^2\int_{0}^{a(t)} \frac{2\gamma}{\gamma-1}\rho^{\gamma-1}(\rho u r^2)_r\mathrm{d}r\\
 \nonumber&&+(1+t)^2\int_{0}^{a(t)}8\pi \rho u\int_{0}^{r}\rho s^{2}\d s\mathrm{d}r\\
 &=&-2(\lambda+2\mu)(1+t)^2\int_{0}^{a(t)}(u_r+\frac{2u}{r})^2 r^2 \mathrm{d}r-(1+t)^2(\rho u^3 r^2+\frac{2}{\gamma-1}\rho^{\gamma}u r^2)|_{r=a(t)}.
\end{eqnarray}
Substituting the above estimates into (\ref{6.2}) yields that
\begin{eqnarray}\label{6.7}
  \nonumber H^\prime(t) &=&-2(\lambda+2\mu)(1+t)^2\int_{0}^{a(t)} (u_r+\frac{2u}{r})^2r^2\mathrm{d}r+\frac{2(2-3(\gamma-1))}{\gamma-1}(1+t)\int_{0}^{a(t)}\rho^{\gamma} r^2\mathrm{d}r \\
  \nonumber&&+2(1+t)\int_{0}^{a(t)} \{4\pi\rho r\int_{0}^{r}\rho s^{2}\d s-\f{4\pi}{r^{2}}(\int_{0}^{r}\rho s^{2}\d s)^2\}\mathrm{d}r\\
 \nonumber&&+ 6(\lambda+2\mu)(1+t) u(a(t),t) a^2(t)-(1+t)^2\f{4\pi u(a(t),t)}{a^2(t)}(\int_{0}^{a(t)}\rho r^{2}\d r)^2\\
 \nonumber &=&-2(\lambda+2\mu)(1+t)^2\int_{0}^{a(t)} (u_r+\frac{2u}{r})^2r^2\mathrm{d}r+\frac{2(5-3\gamma)}{\gamma-1}(1+t)\int_{0}^{a(t)}\rho^{\gamma} r^2\mathrm{d}r \\
\nonumber&&-\f{1}{4\pi}(1+t)\int_{0}^{a(t)}r^2|\Phi_r|^2 \mathrm{d}r+ 6(\lambda+2\mu)(1+t) a^\prime(t) a^2(t)\\
  &&+(1+t)\f{M^2}{4\pi a(t)}-(1+t)^2\f{M^2 }{4\pi}\f{a^\prime(t)}{ a^2(t)},
\end{eqnarray}
where we have used $\Phi_r=-\f{4\pi}{r^2}\int_{0}^{r}\rho s^2 \d s$ and
\begin{eqnarray*}
  4\pi \int_{0}^{a(t)}\rho r\int_{0}^{r}\rho s^{2}\d s\d r&=&2\pi \int_{0}^{a(t)} r^{-1}\partial_r(\int_{0}^{r}\rho s^{2}\d s)^2\d r\\
&=& \f{1}{8\pi}\int_{0}^{a(t)}r^2|\Phi_r|^2 \mathrm{d}r+\f{2\pi }{a(t)}(\int_{0}^{a(t)}\rho r^{2}\d r)^2.
\end{eqnarray*}
By the definition of $H(t)$, we have
\begin{eqnarray}\label{a1}
 \nonumber -\f{1}{4\pi}(1+t)\int_{0}^{a(t)}r^2|\Phi_r|^2 \mathrm{d}r &=&\f{H(t)}{1+t} -\f{1}{1+t}\int_{0}^{a(t)}(r-(1+t) u)^2\rho r^{2}\mathrm{d}r
-\frac{2}{\gamma-1}(1+t) \int_{0}^{a(t)}\rho^{\gamma} r^{2}\mathrm{d}r\\
 &\leq& \f{H(t)}{1+t}-\frac{2}{\gamma-1}(1+t)\int_{0}^{a(t)}\rho^{\gamma} r^{2}\mathrm{d}r.
\end{eqnarray}
Therefore,
\begin{eqnarray}\label{a2}
  \nonumber H^\prime(t) &\leq&\f{H(t)}{1+t}+\frac{2(4-3\gamma)}{\gamma-1}(1+t)\int_{0}^{a(t)}\rho^{\gamma} r^2\mathrm{d}r \\
  \nonumber &&+2(\lambda+2\mu)(1+t) \f{\d}{\d t} a^3(t)+\f{\d}{\d t}(\f{M^2}{4\pi}\f{(1+t)^2}{ a(t)})-\f{M^2 }{4\pi}\f{(1+t)}{ a(t)}.
\end{eqnarray}
Let $Y(t)=H(t)-\f{M^2 }{4\pi}\f{(1+t)^2}{ a(t)}$, then
\begin{equation}\label{a3}
   Y^\prime(t) \leq\f{Y(t)}{1+t}+\frac{2(4-3\gamma)}{\gamma-1}(1+t)\int_{0}^{a(t)}\rho^{\gamma} r^2\mathrm{d}r +2(\lambda+2\mu)(1+t) \f{\d}{\d t} a^3(t).
\end{equation}
For $\f{6}{5}<\gamma\leq \f{4}{3}$ and $M<\overline{M}<M_c$, we obtain the lower bound of $Y(t)$ below by the definition of $H(t)$
and the basic estimate $\eqref{f2.110}$:
\begin{eqnarray}\label{a6}
 \nonumber  Y(t)&=&H(t)-\f{M^2 }{4\pi}\f{(1+t)^2}{ a(t)} \\
 \nonumber &\geq& \frac{2}{\gamma-1}(1+t) ^2\int_{0}^{a(t)}\rho^{\gamma} r^{2}\mathrm{d}r -(1+t) ^2\int_{0}^{a(t)}\f{4\pi}{r^{2}}(\int_{0}^{r}\rho s^{2}\d s)^2\mathrm{d}r-\f{M^2 }{4\pi}\f{(1+t)^2}{ a(t)}\\
 \nonumber &\geq&2(1+t) ^2(\frac{1}{\gamma-1}\int_{0}^{a(t)}\rho^{\gamma} r^{2}\mathrm{d}r-4\pi\int_{0}^{a(t)}\rho r\int_{0}^{r}\rho s^{2}\mathrm{d}s\mathrm{d}r)\\
 &\geq&  (1+t) ^2\frac{1}{\gamma-1}\int_{0}^{a(t)}\rho^{\gamma} r^{2}\mathrm{d}r > 0.
\end{eqnarray}
By \eqref{a3} and \eqref{a6}, we obtain
\begin{equation}\label{aa1}
   Y^\prime(t) \leq 3(3-2\gamma)\f{Y(t)}{1+t} +2(\lambda+2\mu)(1+t) \f{\d}{\d t} a^3(t) ,
\end{equation}
which, together with Gronwall's inequality and \eqref{a6}, implies
\begin{equation}\label{aa2}
 \frac{1}{\gamma-1}\int_{0}^{a(t)}\rho^{\gamma} r^{2}\mathrm{d}r\leq C(1+t)^{7-6\gamma}a^3_1(t),
\end{equation}
since $\f{6}{5}<\gamma\leq\f{4}{3}$ and
\begin{equation}\label{aa3}
  a_1(t):=\max_{s\in[0,t]}a(s)\geq c_0>0.
\end{equation}
Then, we obtain
 \begin{equation}\label{aaa4}
   a_1(t)\geq C(1+t)^{\f{6\gamma-7}{3\gamma}},
\end{equation}
 which deduced from \eqref{aa2} and the fact
\begin{equation}\label{a9}
 \f{M}{4\pi}=\int_{0}^{a_0}\rho_0 r^2\mathrm{d}r= \int_{0}^{a(t)}\rho r^2\mathrm{d}r\leq a(t)^{\frac{3(\gamma-1)}{\gamma}}(\int_{0}^{a(t)}\rho^{\gamma} r^2\mathrm{d}r)^{\frac{1}{\gamma}}.
\end{equation}
By \eqref{a3} and Gronwall's inequality, we obtain
\begin{equation}\label{a4}
 Y(t)\leq (1+t)Y(0)+2(\lambda+2\mu)(1+t)a^3(t)+\frac{2(4-3\gamma)}{\gamma-1}(1+t)\int_{0}^{t}\int_{0}^{a(\tau)}\rho^{\gamma} r^2\mathrm{d}r\d\tau
\end{equation}
For $\gamma = \f{4}{3}$, it is derived from \eqref{a4} that
\begin{equation}\label{a11}
   Y(t)\leq (1+t)Y(0)+2(\lambda+2\mu)(1+t)a^3(t)\leq C(1+t)a^3(t),
\end{equation}
since $a(t)\geq c_0.$
Thus, from \eqref{a3}, we obtain
\begin{equation}\label{a8}
  \int_{0}^{a(t)}\rho^{\gamma}r^2 \d r\leq C (1+t)^{-1}a^3(t).
\end{equation}
The combination of (\ref{a9}) and (\ref{a8}) gives rise to
\begin{equation}\label{6.23}
    a(t)\geq C(1+t)^{\frac{1}{4}},\ \gamma= \frac{4}{3}.
\end{equation}
For $\f{6}{5}<\gamma < \f{4}{3}$, it holds from \eqref{a4} and \eqref{a6} that
\begin{equation}\label{a12}
\frac{1+t}{\gamma-1}\int_{0}^{a(t)}\rho^{\gamma} r^{2}\mathrm{d}r\leq Y(0)+2(\lambda+2\mu)a^3(t)+\frac{2(4-3\gamma)}{\gamma-1}\int_{0}^{t}\int_{0}^{a(\tau)}\rho^{\gamma} r^2\mathrm{d}r\d\tau,
\end{equation}
which, together with Gronwall's inequality, implies
\begin{eqnarray}\label{a13}
 \nonumber\frac{1}{\gamma-1}\int_{0}^{t}\int_{0}^{a(\tau)}\rho^{\gamma} r^2\mathrm{d}r\d\tau\leq&&\frac{Y(0)}{2(4-3\gamma)}[(1+t)^{2(4-3\gamma)}-1]\\
&&+2(\lambda+2\mu)(1+t)^{2(4-3\gamma)}\int_{0}^{t}(1+\tau)^{-9+6\gamma} a^3(\tau)d\tau.
\end{eqnarray}
Claim: for each $\gamma\in (\f{6}{5},\f{4}{3})$ and $\beta\in (\f{2(4-3\gamma)}{3},\f{1}{3\gamma})$, there exists $\{t_n\}$ such that
\begin{equation}\label{a14}
   \f{a(t_n)}{(1+t_n)^\beta}\rightarrow +\infty,~~t_n\rightarrow +\infty.
\end{equation}
Indeed, if \eqref{a14} does not hold, then there is a positive constant $A_\beta$ such that
\begin{equation}\label{a15}
  a(t)\leq A_\beta (1+t)^\beta.
\end{equation}
By \eqref{a9}, \eqref{a12}, \eqref{a13}, \eqref{a15}, we have
\begin{align}\label{a16}
 \nonumber \frac{1}{\gamma-1}(\f{M}{4\pi})^\gamma &\leq Y(0)(1+t)^{2(4-3\gamma)-1}a^{3(\gamma-1)}(t)+2(\lambda+2\mu)(1+t)^{-1}a^{3\gamma}(t)\\
  \nonumber&+4(4-3\gamma)(\lambda+2\mu)(1+t)^{2(4-3\gamma)-1}a^{3(\gamma-1)}(t)\int_{0}^{t}(1+\tau)^{-9+6\gamma} a^3(\tau)d\tau\\
  \nonumber&\leq Y(0) A^{3(\gamma-1)}_\beta (1+t)^{2(4-3\gamma)-1+3\beta(\gamma-1)}+2(\lambda+2\mu)A^{3\gamma}_\beta(1+t)^{3\beta\gamma-1}\\
 \nonumber&+ 4(4-3\gamma)(\lambda+2\mu)A^{3\gamma}_\beta(1+t)^{2(4-3\gamma)-1+3\beta(\gamma-1)}\int_{0}^{t}(1+\tau)^{-9+6\gamma+3\beta} d\tau\\
  \nonumber&\leq Y(0) A^{3(\gamma-1)}_\beta(1+t)^{2(4-3\gamma)-1+3\beta(\gamma-1)}+2(\lambda+2\mu)A^{3\gamma}_\beta(1+t)^{3\beta\gamma-1}\\
&+4(4-3\gamma)(\lambda+2\mu)A^{3\gamma}_\beta\f{(1+t)^{3\beta\gamma-1}}{3\beta-2(4-3\gamma)}\rightarrow 0,~ as~~t\rightarrow+\infty.
\end{align}
This is a contradiction. Thus, the time expanding rate \eqref{a14} of the free boundary holds.

\bigskip

\noindent {\bf Acknowledgements:}
Thanks the referees for their insightful comments
and suggestions to improve this paper. The research is supported by the NNSFC (grants No.
11671384 and 11225102), NSFC-RGC Grant 11461161007 and by the Key Project of Beijing Municipal Education
Commission No. CIT$\&$TCD20140323.


\medskip












\addcontentsline{toc}{section}{References}
\begin{thebibliography}{99}

\bibitem{C1938}S. Chandrasekhar, An Introduction to the Study of Stellar Structures. University of Chicago Press, Chicago, 1938.

\bibitem{DLYY}
Y. Deng, T. P. Liu, T. Yang, Z. A. Yao, Solutions of Euler-Poisson equations for
gaseous stars. \textit{Arch. Ration. Mech. Anal.} \textbf{164} (2002), 261-285.

\bibitem{DL2015}
Q. Duan, H. L. Li, Global existence of weak solution for the compressible
Navier-Stokes-Poisson system for gaseous stars. \textit{J. Diff. Eqs. }\textbf{259} (2015), 5302-5330.

\bibitem{Ev1998}
L. C. Evans, Partial differential equations. Graduate Studies in
Mathematics. American Mathematical Society, Providence, RI,
1998.

\bibitem{SG}
S. Gao, Global solutions to the Navier-Stokes-Poisson equations for self-gravitating gaseous stars.
Thesis (Ph.D.) -Northwestern University. 2010. 89 pp.

\bibitem{GLX2012}
Z. H. Guo, H. L. Li, Z. P. Xin, Lagrange structure and dynamics for spherically symmetric compressible Navier-
Stokes equations. \textit{Comm. Math. Phys.}
\textbf{309}(2) (2012), 371-412.


\bibitem{JJ2008}
J. Jang, Nonlinear instability in gravitational Euler-Poisson systems for $\gamma=\f{6}{5}$.
\textit{Arch. Rational Mech. Anal.} \textbf{ 188} (2008), 265-307.

\bibitem{JJ2010}
J. Jang, Local well-posedness of dynamics of viscous gaseous stars.
\textit{Arch. Rational Mech. Anal.} \textbf{ 195} (2010), 797-863.

\bibitem{JT2013}
J. Jang, I. Tice, Instability theory of the Navier-Stokes-Poisson equations.
\textit{Anal. PDE} \textbf{ 5} (2013), 1121-1181.


\bibitem{JZ2001}
S. Jiang, P. Zhang, On spherically symmetric solutions of the compressible isentropic
Navier-Stokes equations. \textit{Commun. Math. Phys.} \textbf{215} (2001), 549-581.

\bibitem{KLL} H. H. Kong, H. L. Li, C. C. Liang, Global solutions to 3D isentropic compressible Navier-Stokes equations with free boundary. \textit{Bulletin of the Institute of Mathematics Academia Sinica (New Series)}
Vol. 10 (2015), No. 4, pp. 575-613.

\bibitem{L1997}
S. S. Lin, Stability of gaseous stars in spherically symmetric motions.
\textit{SIAM J. Math. Anal.} \textbf{28} (1997), 539-569.


\bibitem{LXY2000}T. Luo, Z. P. Xin, T. Yang, Interface behavior of compressible Navier-Stokes equations with vacuum. \textit{SIAM J. Math. Anal. } \textbf{31} (2000), 1175-1191.

\bibitem{LXZ}
T. Luo, Z. P. Xin, H. H. Zeng, Nonlinear asymptotic stability of the Lane-Emden solutions for the viscous gaseous star problem with degenerate density dependent viscosities. arXiv:1507.01069.

\bibitem{LXZ1} T. Luo, Z. P. Xin, H. H. Zeng, Nonlinear asymptotic stability of the Lane-Emden solutions for the viscous gaseous star problem, arXiv:1506.03906.
 \bibitem{ZF2009} T. Zhang, D. Y. Fang, Global behavior of spherically symmetric Navier-Stokes-Poisson system with degenerate viscosity coefficients. \textit{Arch. Ration. Mech. Anal.} \textbf{191} (2009), no. 2, 195-243.

\end{thebibliography}
\end{document}